# Distributionally robust decision-making under ambiguity: case study of water environmental management

Hidekazu Yoshioka[a, *], Motoh Tsujimura[b], Yumi Yoshioka[c], Ayumi Hashiguchi[d]

[a] Graduate School of Advanced Science and Technology, Japan Advanced Institute of Science and Technology, Japan

[b] Graduate School of Commerce, Doshisha University, Japan

[c] Faculty of Applied Biological Sciences, Gifu University, Japan

[d] Faculty of Environmental, Life, Natural Science and Technology, Okayama University, Japan

* Corresponding author: Associate Professor, Graduate School of Advanced Science and Technology, Japan Advanced Institute of Science and Technology, Asahidai 1-1, Nomi, Japan. Tel: +81-761-51-1745 , E-mail: yoshih@jaist.ac.jp

**Abstract**

Decision-making under uncertainty is ubiquitous in environmental project planning. Environmental processes such as a streamflow discharge often present a subexponential memory, where the autocorrelation persists for a long time. In addition, optimization problems driven by environmental processes encounter the issue of model ambiguity because of a lack of sufficient data for model identification. To facilitate decision-making for the management of aquatic environments (e.g., flood mitigation, water abstraction for hydropower generation), we formulate a unified distributionally robust stochastic optimization problem based on a mixed moving average (MMA) process. The MMA process is a superposition of infinite-dimensional affine stochastic processes that is seemingly complex, but the affine property helps with the formulation and computation of the optimization. Our problem is based on a convex objective with a nonsmooth conditional value-at-risk measure. We present a convergent regularization to obtain its smooth and strictly convex counterpart. The model ambiguity is represented as a distortion of the probability density of the target dynamics, and it is penalized by a divergence with which the optimization problem remains convex and becomes computable. As a case study, we apply the optimization problem to two cases with identified parameter values. The performance of the optimized dynamics is evaluated through a statistical simulation. This paper serves as a multidisciplinary work covering both the theory and application of distributionally robust optimization.

**Statements and declarations**

***Acknowledgments***
This work was supported by Japan Society for the Promotion of Science (22K14441, 22H02456), Environmental Research Projects from the Sumitomo Foundation (203160), and Kurita Water and Environment Foundation (21K018).

***Availability of data and material***
Data will be available upon a reasonable request to the corresponding author.

***Competing interests***
The authors have no competing interests.

***Author contributions***
Hidekazu Yoshioka: Conceptualization, Methodology, Software, Formal analysis, Data Curation, Writing – original draft preparation, Writing – review and editing, Supervision, Project administration, Funding acquisition
Motoh Tsujimura: Formal analysis, Writing – review and editing
Yumi Yoshioka: Investigation, Writing – review and editing, Visualization
Ayumi Hashiguchi: Investigation, Writing – review and editing

## 1. Introduction

### 1.1 Background

A crucial step for sustainable development is coexistence with the environment and ecosystem. Decision-making under uncertainty is ubiquitous in environmental project planning because of the need to optimize inherently stochastic processes. State-space models have been developed to realize sustainable crop production in greenhouses under uncertain weather conditions (Chen and You, 2021; Hu and You, 2022). Management of marine fisheries requires the consideration of diverse factors, such as economic and social development, governance and policy, biology and the environment, and logistics (Akbari et al., 2022). Land use management has been analyzed as a principal factor for species conservation and water pollution (Aggarwal et al., 2022; Epanchin-Niell et al., 2022). Rivers are the core of the water–energy–ecosystem nexus, and they are among the most severely threatened environments because of anthropogenic activities (Yin et al., 2022; Zhang et al., 2022a). For example, water diversion improves water supply resources while degrading regional water quality (Zhang et al., 2022b). Thus, river management has been a long-standing engineering issue.

Stochastic optimization provides a flexible and rigorous mathematical framework for decision-making under uncertainty (Shapiro et al., 2021) so that target dynamics can be managed. In a stochastic optimization problem, the decision variables are optimized so that an objective based on some expectation is minimized subject to the constraints of the dynamics and the variables themselves. Zhang and Xi (2020) proposed an interval two-stage optimization problem for water resource management. Salgado et al. (2022) discussed flood control through dam–reservoir systems for large hydropower generation schemes, such as the Three Gorges Reservoir in China. Lazzaro and Botter (2015) developed a statistical model to evaluate the performance of smaller-scale run-of-river hydropower generation projects. Huuki et al. (2022) used nonlinear programming to study the ecologically friendly operation of smaller hydropower projects.

In practice, optimization problems driven by environmental processes, which often include hydrological variables such as streamflow discharge, are difficult to describe completely because of the lack of sufficient data to determine target dynamics (Bacci et al., 2022; Hah et al., 2022; Engman et al., 2021). In such cases, the problem is considered ambiguous. Distributionally robust stochastic optimization is an alternative to conventional optimization (Postek et al., 2016), where the target dynamics can be distorted in the sense of probability distribution. In this framework, the model ambiguity is represented as the potential difference between the benchmark and distorted (worse) models, which is measured as a divergence (Lin et al., 2022) or distance (Gao and Kleywegt, 2022). The distributionally robust nature leads to a game problem instead of a minimization problem, which introduces ambiguity-induced nonlinearity even if the problem is originally linear. Nevertheless, this approach is widely used in engineering studies owing to its high computational tractability and well-posedness under certain assumptions. Examples of its application include risk-averse reservoir operation to reduce flood risk (Gauvin et al., 2017), robust optimization of hydro–wind–solar energy systems (Jin et al., 2022), Q-learning of the Markov decision process (Bäuerle and Glauner, 2021), and optimization of the uncertain fractional partial differential equation (Anti et al., 2021). Thus, distributionally robust optimization is a pivotal tool for decision-making in various research fields.

Another critical issue for optimization problems is mathematical modeling of system dynamics. Actual environmental processes often exhibit long (subexponential) memory that persists significantly longer than processes with exponential memory do. Typical examples include streamflow discharge (Dey and Mujumdar, 2022; Yoshioka et al., 2022a), water quality indices (Rozental and Tambieva, 2021), and pollutant and greenhouse gas concentrations in the air (Bukhari et al., 2022; Proietti and Maddanu, 2022). Candidate models for describing system dynamics exhibiting subexponential memory include time-fractional stochastic processes and the superposition of processes possessing multiple timescales (Barndorff-Nielsen et al., 2018). However, the application of such models to stochastic optimization is still rare despite their relevance to engineering problems.

### 1.2 Objective and contribution

The objective of this paper is to present a unified framework of distributionally robust optimization for continuous-time stochastic processes with subexponential memory. Our focus is on its application to environmental planning in river environments as a case study. We consider that it is also applicable to planning problems in other research fields owing to its mathematical universality by suitably interpreting the meaning of stochastic processes. Our approach is based on previous studies (Yoshioka et al., 2023b-c) that discussed static optimization problems. The differences between the present study and previous studies are also explained below.

We target a mixed moving average (MMA) process, which is formally a superposition of infinitely many mutually independent stochastic differential equations (SDEs) with different timescales (Rajput and Rosinski, 1989). MMA processes generate subexponential memory by spectrally mixing exponential memories. This is fundamentally different from fractional models (e.g., Volterra), in which an algebraic memory kernel plays a vital role (Hainaut, 2022; Shi et al., 2022). The simplest MMA process is the superposition of the Ornstein–Uhlenbeck (supOU) type, where Ornstein–Uhlenbeck (OU) processes are integrated with respect to the reversion speed in a phase space through Lévy bases (Barndorff-Nielsen, 2001). Its time-discrete counterpart is the random-coefficient autoregressive model (Pilipauskaitė et al., 2020). The superposition of continuous-state branching processes (called the supCBI process) as an affine generalization of the supOU process was previously proposed (Yoshioka, 2022) to describe clustered flood events (e.g., Fenicia and McDonnell, 2022). MMA processes, especially those based on affine processes such as OU and CBI processes, are attractive because they handle both exponential and subexponential memories relevant for many time series data in the problems of finance and economics in a unified manner. In addition, MMA processes have statistical moments in the closed form. The stochastic process model here is more general than the supOU process used in the previous study (Yoshioka et al., 2023c). In this paper, we use supCBI processes because of their flexibility and tractability. The characteristic functions of these processes have computable forms. Owing to this useful property, their probability density functions (PDFs) become accessible without losing the tractability of the model.

**Figure 1** shows the flow of our optimization modeling. We formulate a unified mathematical model of an optimization problem to cover two cases: high and low flows. The risk of diverting too much or too little streamflow is evaluated by a conditional value-at-risk (CVaR) measure that harmonizes with the convex optimization because of the dual representation (Ang et al., 2018). We employ CVaR because of its wide use in related studies such as smart energy management (Zhai et al., 2022), natural resource management (Mir et al., 2022), and insurance risk management (Ghossoub et al., 2022). A CvaR term has also been used in the optimization problem of Yoshioka et al. (2023b), while our problem uses a pair of CVaRs and hence is slightly generalized.

Model ambiguity is evaluated by using a Radon–Nikodým derivative between the benchmark and distorted (worse) models (Liu et al., 2022). The ambiguity is then penalized by a divergence, which keeps the optimization problem convex. The form of the optimization problem presented in this paper is a variant of those in the previous studies; a study that investigated a problem without uncertainties and hence without distributionally-robust nature (Yoshioka et al., 2023b), and a study that analyzed distributionally-robust optimization computationally (Yoshioka et al., 2023c); this paper presents several theoretical results along with applications. More specifically, this paper discusses conditions to well-pose the optimization problem.

We present a convergent regularization method to smoothen the CVaR term, which yields a smooth and strictly convex optimization problem. The strict convexity suggests that our problem yields a globally optimal decision variable under ambiguity that is computable by existing methods such as gradient descent (e.g., Iyengar and Ma, 2013), the accelerated inertial method (e.g., Wibisono et al., 2016), or its generalization (Liang et al., 2020). We use the inertial method because it is more efficient than gradient descent (Calder and Yezzi, 2019; Dambrine et al., 2022; Duruisseaux and Leok, 2022; He et al., 2022; Wang et al., 2022). The convergence of the optimization problems is analyzed as well.

We apply the optimization problem to two case studies where the streamflow is diverted. Parameter values of the supCBI process were identified from actual data, and the optimization problem was discretized by using the Markovian lift (Yoshioka et al., 2022a; Yoshioka et al., 2023a) to reduce an MMA process to a finite-dimensional system of SDEs. We found that a nonproportional diversion rule such as that in Razurel et al. (2016) and subsequently applied to several studies (Razurel et al., 2018; Perona et al., 2021) emerges in our solution, especially for streamflow diversion to generate hydropower. Thus, our contribution suggests the rationality of the existing parameterized scheme. The performance of the optimization result is analyzed by a Monte Carlo simulation. Consequently, this paper serves as a multidisciplinary work covering both the theory and application of distributionally robust optimization in environmental and operations research areas. In particular, the modeling and computational methods used to address MMA processes are transferable to problems in research areas other than environmental management, as these methods are based on versatile mathematical descriptions.

The rest of this paper is organized as follows. **Section 2** presents the MMA processes and their basic results. **Section 3** formulates the optimization problem and analyzes its well-posedness and convergence in different cases. **Section 4** presents applications of our optimization framework. **Section 5** summarizes our results and presents future perspectives. The **Appendices** present supplementary results, including proofs of propositions and supporting data.

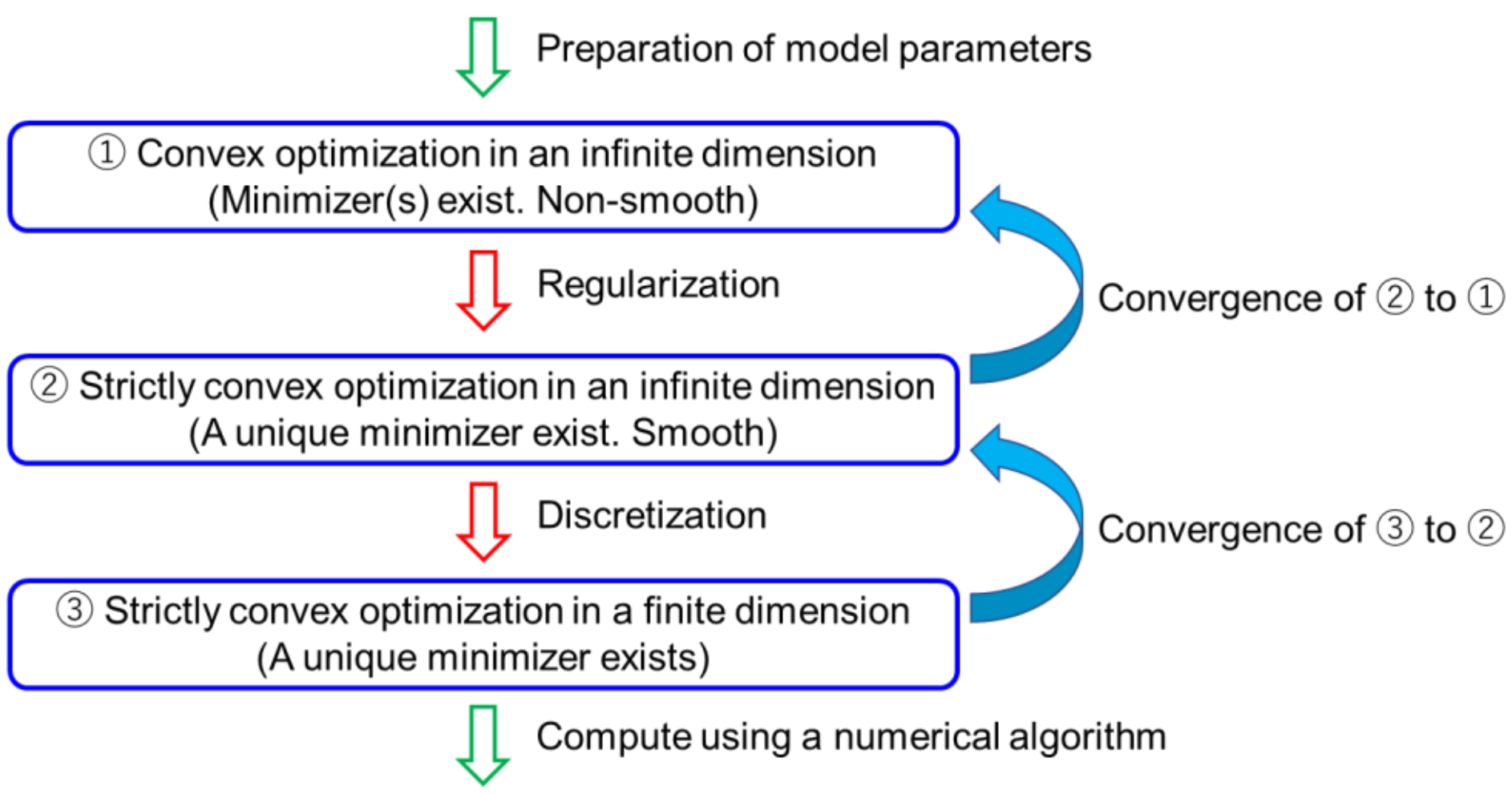

**Figure 1.** The flow of our optimization modeling.

## 2. Mathematical model

This section introduces mathematical models for streamflow discharge and summarizes their basic results.

### 2.1 Stationary MMA processes

The MMA process is formally a superposition, namely, an integration with respect to Lévy bases, of mutually independent SDEs shearing a generic form. The convergence of superposition is justified in the sense of distribution (Yoshioka, 2022), where the MMA process is approximated as a summation of a finitely many mutually independent SDEs. This technique is called Markovian lift and has been used to reformulate a non-Markovian process into a system of Markovian processes with larger dimensions but easier to handle (Abi Jaber, 2019; Siegle et al., 2010); moreover, the MMA process admits explicit moments and autocorrelation functions that can be implemented efficiently in applied problems. The explanation of the MMA processes below is based on the methods of Barndorff-Nielsen (2001), Barndorff-Nielsen et al. (2011), and Yoshioka (2022). We obtain several estimates to be used in **Section 3**.

We assume a complete probability space as usual in conventional stochastic calculus. We consider supOU and supCBI processes as tractable MMA processes, both of which generate paths with successive jumps decaying subexponentially. We use the supCBI process in the application but the explanation in this section starts from the supOU process as a special case to better understand MMA processes.

#### 2.1.1 supOU process

The supOU process is an integration of infinitely many OU processes through Lévy bases as a model of Lévy random fields (**Definition 1** below). Here, $\mathbb{B}\left(\left(0,+\infty\right)\times\mathbb{R}\right)$ is the Borel sigma algebra on the product space $(0,+\infty)\times\mathbb{R}$, i is the imaginary unit ($\mathrm{i}^2=-1$), and $\mathbb{E}$ is the expectation.

***Definition 1*** *A family of real scalar random variables* $\Lambda=\{\Lambda(B): B\in\mathbb{B}((0,+\infty)\times\mathbb{R})\}$ *satisfying the conditions (a)-(e) below is called Lévy bases.*

*(a) The distribution of* $\Lambda(B)$ *is infinitely divisible for any* $B\in\mathbb{B}((0,+\infty)\times\mathbb{R})$.

*(b) For any* $n\in\mathbb{N}$ *and any mutually disjoint sets* $B_1,B_2,...,B_n\in\mathbb{B}((0,+\infty)\times\mathbb{R})$, *random variables* $\Lambda(B_1),\Lambda(B_2),...,\Lambda(B_n)$ *are mutually independent.*

*(c) For any arbitrary mutually disjoint sets* $\{B_i\}_{i=1,2,3,...}$ *satisfying* $B_i\in\mathbb{B}((0,+\infty)\times\mathbb{R})$ *(* $i\in\mathbb{N}$ *) and* $\bigcup_{i\in\mathbb{N}}B_i\in\mathbb{B}((0,+\infty)\times\mathbb{R})$, $\sum_{i=1}^{+\infty}\Lambda(B_i)$ *converges almost surely (a.s.) and* $\Lambda\left(\bigcup_{i\in\mathbb{N}}B_i\right)=\sum_{i=1}^{+\infty}\Lambda(B_i)$.

*(d)* $\mathbb{E}[\exp(\mathrm{i}u\Lambda(B))]=\exp(\varphi(u)\Pi(B))$ *for any* $u\in\mathbb{R}$ *and* $B\in\mathbb{B}((0,+\infty)\times\mathbb{R})$, *where* $\Pi$ *is the product measure of a probability measure* $\pi$ *and a Lebesgue measure, and* $\varphi$ *is a characteristic function of a non-decreasing pure-jump Lévy process having bounded variations.*

*(e) The probability measure* $\pi$ *satisfies the inverse-moment boundedness* $\int_0^{+\infty}\rho^{-1}\pi(\mathrm{d}\rho)<+\infty$.

The conditions (a)-(c) are common in conventional Lévy bases. The condition (d) limits Lévy bases to time-homogenous ones like Lévy processes, where the assumption on the characteristic function turns out to be satisfied in our applications. This condition is necessary for the stationarity of MMA processes, physically meaning that the distributed reversion parameter is not significantly accumulated near the origin corresponding to the infinitely fast reversion. The condition (e) is included in this paper to obtain a MMA process with bounded statistical moments (e.g., Eq. (7)).

We introduce the supOU process $X=(X_t)_{t\in\mathbb{R}}$ as a stationary scalar process, which is the streamflow discharge (the volume of water passing through a river cross-section):

$$X_t=\underline{X}+\int_0^{+\infty}\int_{-\infty}^{t}\exp(-\rho(t-s))\Lambda(\mathrm{d}\rho,\mathrm{d}s),\quad t\in\mathbb{R} \tag{1}$$

with a prescribed constant $\underline{X}\geq 0$ called the minimum discharge (because the second term in (1) is nonnegative), where the outer integration is for the reversion parameter $\rho>0$, whereas the inner integration is for time $s$. The second term in (1) represents jump-driven flood events in which each sudden discharge increase is followed by gradual recession. The process (1) is understood as a formal integration of the measure-valued mutually independent OU processes

$$Y_t(\mathrm{d}\rho)=\int_{-\infty}^{t}\exp(-\rho(t-s))\Lambda(\mathrm{d}\rho,\mathrm{d}s),\quad t\in\mathbb{R}, \tag{2}$$

leading to

$$X_t = \underline{X} + \int_0^{+\infty} Y_t(\mathrm{d}\rho),\ \ t \in \mathbb{R}. \tag{3}$$

It is an integration of mutually independent OU processes with different reversion parameter values. The coexistence of different reversion speeds indicates the existence of multiple time scales in the physical processes that generate streamflow discharge and appear as subexponential memory (Yoshioka et al., 2023a).

We can rewrite (1) as follows:

$$X_t = \underline{X} + \int_0^{+\infty}\int_0^{+\infty}\int_{-\infty}^{t} \exp(-\rho(t-s)) z\mu(\mathrm{d}z, \mathrm{d}\rho, \mathrm{d}s),\ \ t \in \mathbb{R} \tag{4}$$

using a Poisson random measure on $(0,+\infty)\times(0,+\infty)\times\mathbb{R}$ with the intensity $Av(\mathrm{d}z)\pi(\mathrm{d}\rho)\mathrm{d}s$ and a constant $A>0$. Here, $v$ is a Lévy measure of a nondecreasing pure-jump Lévy process with bounded variations, and the outset integration of (4) is for the jump size $z>0$.

The advantage of the supOU process is its tractability in that the characteristic and moment-generating functions are found analytically, from which a variety of statistics can be obtained. The probability density of a random variable is obtained through a Fourier inversion of the corresponding characteristic function (e.g., Chapter 5 of Hainaut (2022)). Hence, the remarkable tractability of the supOU process allows for efficient computation of its PDF. Furthermore, the autocorrelation function with a subexponential decay of the process is available in a closed form (e.g., **Appendix B**).

The characteristic function of a random variable $X$ is denoted as $c_X(\xi)$, i.e.,

$$c_X(\xi) = \mathbb{E}\left[e^{\mathrm{i}\xi X}\right],\ \ \xi \in \mathbb{R}, \tag{5}$$

which always exists because $\left|e^{\mathrm{i}\xi X}\right| \leq 1$, while the moment-generating function

$$M_X(\xi) = \mathbb{E}\left[e^{\xi X}\right],\ \ \xi \in \mathbb{R} \tag{6}$$

may not occur if $\xi$ is large. The characteristic function of the supOU process at a stationary state is as follows (e.g., Barndorff-Nielsen, 2001):

$$c_{X_t}(\xi) = e^{\mathrm{i}\xi\underline{X}} \exp\left( A\int_0^{+\infty} \frac{\pi(\mathrm{d}\rho)}{\rho} \int_0^{+\infty}\int_0^{+\infty} \left(\exp\left(\mathrm{i}\xi z e^{-t}\right) - 1\right) v(\mathrm{d}z)\mathrm{d}t \right),\ \ \xi \in \mathbb{R}, \tag{7}$$

while the moment-generating function, if it exists, is

$$M_{X_t}(\xi) = e^{\xi\underline{X}} \exp\left( A\int_0^{+\infty} \frac{\pi(\mathrm{d}\rho)}{\rho} \int_0^{+\infty}\int_0^{+\infty} \left(\exp\left(\xi z e^{-t}\right) - 1\right) v(\mathrm{d}z)\mathrm{d}t \right),\ \ \xi \in \mathbb{R}. \tag{8}$$

The moment-generating function plays a key role in this paper because it enables us to prove the existence of exponential moments of MMA processes. More specifically, the boundedness of exponential moments becomes essential for the well-posedness and convergence of our optimization problem. We can further proceed in the tempered stable case with $v(\mathrm{d}z) = z^{-(\alpha_v+1)}\exp(-\beta_v z)\mathrm{d}z$ ($\alpha_v < 1$, $\beta_v > 0$) (e.g., Yoshioka, 2022) used in our application:

$$\left(A\int_0^{+\infty}\frac{\pi(\mathrm{d}\rho)}{\rho}\right)^{-1}\ln M_{X_t}(\xi)=\int_0^{+\infty}\int_0^{+\infty}\left(\exp\left(\xi z e^{-t}\right)-1\right)v(\mathrm{d}z)\mathrm{d}t$$
$$=\frac{\Gamma(1-\alpha_v)}{\alpha_v}(\beta_v)^{\alpha_v}\int_0^{+\infty}\left(1-\left(1-\frac{\xi e^{-t}}{\beta_v}\right)^{\alpha_v}\right)\mathrm{d}t, \quad \xi\in\mathbb{R}. \tag{9}$$

The last integral is bounded if $\xi\le\beta_v$, and it diverges to $+\infty$ otherwise. This implies that $\xi$ must be small to ensure the existence of $M_{X_t}(\xi)$. Consequently, we obtain **Proposition 1** below.

***Proposition 1*** *For the supOU process with the tempered stable* $v(\mathrm{d}z)$*, if* $\xi\le\beta_v$ *then* $M_{X_t}(\xi)$ *exists and is given by* (8).

We also obtain **Proposition 2** concerning other exponential moments.

***Proposition 2*** *For the supOU process with the tempered stable* $v(\mathrm{d}z)$*, if* $\xi\le\beta_v$ *then*

$$\mathbb{E}\left[X_t e^{\xi X_t}\right]=A' M_{X_t}(\xi)\frac{1}{\xi}\left(1-\left(1-\frac{\xi}{\beta_v}\right)^{\alpha_v}\right)<+\infty \tag{10}$$

*and*

$$\mathbb{E}\left[X_t^2 e^{\xi X_t}\right]=A'\frac{\partial}{\partial\xi}\left[M_{X_t}(\xi)\frac{1}{\xi}\left(1-\left(1-\frac{\xi}{\beta_v}\right)^{\alpha_v}\right)\right]<+\infty \tag{11}$$

*with* $A'=A\int_0^{+\infty}\frac{\pi(\mathrm{d}\rho)}{\rho}\frac{\Gamma(1-\alpha_v)}{\alpha_v}(\beta_v)^{\alpha_v}$.

### 2.1.2 supCBI process

The supCBI process is a self-exciting generalization of the supOU process recently discussed in Yoshioka (2023). It is understood through a Markovian lift as an approximation that converges in the sense of distribution. This process at a stationary state has the characteristic function

$$c_{X_t}(\xi)=e^{\mathrm{i}\xi\underline{X}}\exp\left(A\int_0^{+\infty}\frac{\pi(\mathrm{d}\rho)}{\rho}\int_0^{+\infty}\int_0^{+\infty}\left(\exp\left(\mathrm{i}\varpi(t)z\right)-1\right)v(\mathrm{d}z)\mathrm{d}t\right),\quad \xi\in\mathbb{R}. \tag{12}$$

Here, $\varpi(t)$ is the unique smooth solution, if it exists, to the initial value problem

$$\mathrm{i}\frac{\mathrm{d}\varpi(t)}{\mathrm{d}t}=-\mathrm{i}\rho\varpi(t)+B\int_0^{+\infty}\left(\exp\left(\mathrm{i}\varpi(t)z\right)-1\right)v(\mathrm{d}z),\quad t>0,\quad \varpi(0)=\xi. \tag{13}$$

We interpret the supCBI process as a weak limit (i.e., in the sense of distribution) of the finite-dimensional case

$$X_t=\underline{X}+\sum_{i=1}^{n}Y_t(\rho_i),\quad t\in\mathbb{R} \tag{14}$$

under $n \to +\infty$ with each $Y_t(\rho_i)$ ($1 \le i \le n$) being the unique strong solution to

$$\mathrm{d}Y_t(\rho_i) = -\rho_i Y_t(\rho_i)\mathrm{d}t + \int_0^{c_i A + r_i B Y_t(\rho_i)} \int_0^{\infty} z \mu_i(\mathrm{d}z, \mathrm{d}u, \mathrm{d}s), \quad t \in \mathbb{R} \tag{15}$$

such that each $Y_t(\rho_i)$ and $Y_t(\rho_j)$ are independent of each other when $i \ne j$. Here, $B \ge 0$ is a positive parameter, $\{\mu_i\}_{1 \le i \le n}$ is a sequence of mutually independent Poisson random measures on $(0,+\infty) \times (0,+\infty) \times \mathbb{R}$ with the compensated measure $v(\mathrm{d}z)\mathrm{d}u\mathrm{d}s$, $\{\rho_i\}_{1 \le i \le n}$ is a positive and strictly increasing sequence of reversion parameters, and the sequence $\{c_i\}_{1 \le i \le n}$ is a discrete counterpart of the probability density $\pi$ and defines the discrete probability measure $\sum_{i=1}^{n} c_i \delta_{(\rho=\rho_i)}$, with $\delta_{(\rho=\rho_i)}$ being the Dirac delta at $\rho = \rho_i$. A sufficient condition for the convergence of (14), the system corresponding to the Markovian lift, has been obtained (Yoshioka et al., 2023a). The parameter $B$ modulates the strength of the self-excitedness such that a larger $B$ implies stronger self-excitation and disappears if $B = 0$ at which the supCBI process is reduced to a supOU process. The parameter $B$ should satisfy $BM_1 < 1$ ($M_1 = \int_0^{+\infty} z v(\mathrm{d}z)$) for stationarity (Yoshioka, 2022).

***Remark 1*** Our Markovian lift differs from the existing ones for SDEs with distributed delay (Abi Jaber, 2019; Siegle et al., 2010). Our supCBI process discretizes the randomized reversion parameter in the Lévy bases and gives a mutually independent system of SDEs. In contrast, those in the literature approximate a memory kernel with a nonexponential delay by a mixture of exponential kernels. As a result, the resulting Markovian system has SDEs driven by common noise.

The existence of a moment-generating function of the supCBI process holds true for small $\xi$. The same applies to the moment of the form $\mathbb{E}\left[X_t e^{\xi X_t}\right]$ as shown in **Proposition 3** below.

***Proposition 3*** *For the supCBI process with the tempered stable* $v(\mathrm{d}z)$, *if* $\xi \le \beta_c \equiv \beta_v \left(1 - (BM_1)^{\frac{1}{1-\alpha_v}}\right) \in (0, \beta_v)$ *and* $B$ *is sufficiently small, then*

$$\mathbb{E}\left[e^{\xi X_t}\right] = e^{u\underline{X}} \exp\left( A \int_0^{+\infty} \frac{\pi(\mathrm{d}\rho)}{\rho} \int_0^{+\infty} \int_0^{+\infty} \left(\exp(\varpi(t) z) - 1\right) v(\mathrm{d}z) \mathrm{d}t \right) \tag{16}$$

*exists with* $\varpi(t)$ ($t \ge 0$) *being a smooth curve uniquely determined from the initial value problem:*

$$\frac{\mathrm{d}\varpi(t)}{\mathrm{d}t} = -\varpi(t) + B \int_0^{+\infty} \left(\exp(\varpi(t) z) - 1\right) v(\mathrm{d}z), \quad \varpi(0) = \xi, \quad t > 0. \tag{17}$$

***Remark 2*** The quantity $\mathbb{E}\left[X_t^k e^{\xi X_t}\right]$ ($k = 1, 2$) also exists if $\xi \le \beta_c$ and $B$ issmall.

## 3. Optimization problem

We formulate and analyze stochastic optimization problems of MMA processes under ambiguity. Although our focus is on water resource management, the framework presented below can be applied to problems in other research areas because the mathematical results in this section do not strongly rely on the specific problem background.

### 3.1 Benchmark case

Given a measurable, nonnegative, bounded, continuous PDF $p$ of the target dynamics, we first formulate an optimization problem without model ambiguity. The expectation of a function $f = f(X)$ of the nonnegative scalar random variable $X$ is denoted as $\mathbb{E}\left[f(X)\right] = \int f(x) p(x) \mathrm{d}x$, where the domain of integration is indicated when necessary. The upper- and lower-CVaR measures for the generic nonnegative integrable random variable $Z \ge 0$ are defined as follows ($\alpha, \beta \in (0,1)$ are quantile levels):

$$\mathrm{CVaR}_\alpha(Z) = \frac{\mathbb{E}\left[Z\mathbb{I}(Z \le Z_\alpha)\right]}{\mathbb{E}\left[\mathbb{I}(Z \le Z_\alpha)\right]} \text{ and } \mathrm{CVaR}^{1-\beta}(Z) = \frac{\mathbb{E}\left[Z\mathbb{I}(Z \ge Z_\beta)\right]}{\mathbb{E}\left[\mathbb{I}(Z \ge Z_\beta)\right]}. \tag{18}$$

Here, $Z_\gamma$ is the quantile value corresponding to the quantile level $\gamma \in (0,1)$ being uniquely determined by $\gamma = \mathbb{E}\left[\mathbb{I}(Z \le Z_\gamma)\right]$, and $\mathbb{I}$ is the indicator function (e.g., $\mathbb{I}(Z \le Z_\alpha) = 1$ if $Z \le Z_\alpha$ and $\mathbb{I}(Z \le Z_\alpha) = 0$ otherwise). The upper- and lower-CVaR measures are conditional expectations to evaluate the lower- and upper-tail risks of the random variable. By the dual representation formula (e.g., Section 3.4 in Ang et al., 2018), we have the application-oriented forms

$$\mathrm{CVaR}_\alpha(Z) = -\inf_{u \ge 0}\left\{-u + \frac{1}{\alpha}\mathbb{E}\left[\max\{u - Z, 0\}\right]\right\} \tag{19}$$

and

$$\mathrm{CVaR}^{1-\beta}(Z) = \inf_{w \ge 0}\left\{w + \frac{1}{1-\beta}\mathbb{E}\left[\max\{Z - w, 0\}\right]\right\}. \tag{20}$$

We assume that $p$ is positive (a PDF is possibly nonnegative) in $\Xi = (\underline{X}, +\infty)$. The decision variable is the diversion ratio $c = c(x)$ as a bounded measurable function of $x$. Physically, the discharge $x$ from the upstream channel is diverted or taken with the amount $cx$ and the remaining $(1-c)x$ flows down the main channel. The admissible set of decision variables is chosen as follows:

$$\mathfrak{C} = \left\{c \middle| c : \Xi \to [0,1],\ c \in L_p^2(\Xi)\right\}, \tag{21}$$

where $L_p^2(\Xi)$ is the space of functions $c : \Xi \to [0,1]$ such that

$$\int c^2(x)p(x)\mathrm{d}x\left(\leq\int p(x)\mathrm{d}x=1\right)<+\infty, \tag{22}$$

and is equipped with the natural norm $\|c\|=\sqrt{\int c^2(x)p(x)\mathrm{d}x}$. The space $L_p^2(\Xi)$ is a Banach space (e.g., Theorem 1.3. of Kufner and Opic (1984)) and, further, a Hilbert space by square integrability. The space $\mathfrak{C}$ is a convex, closed, and bounded subspace of $L_p^2(\Xi)$. Indeed, the convexity and boundedness are trivial. For closedness, applying the dominated convergence theorem (e.g., Theorem 4.2 of Brezis (2011)) to an arbitrary sequence $c_n\in\mathfrak{C}$ ($n=1,2,3,...$) (uniformly bounded by definition) converging pointwise to some $c:\Xi\to[0,1]$ almost everywhere in $\Xi$ yields $c\in\mathfrak{C}$ because

$$\lim_{n\to+\infty}\int c_n^2(x)p(x)\mathrm{d}x=\int c^2(x)p(x)\mathrm{d}x\leq 1<+\infty. \tag{23}$$

The benchmark optimization problem aims at minimizing an objective to obtain a compromising policy among the diversion effort and the risk of too much/less diverting the flow. Set the objective

$$J(c)=\mathbb{E}\left[g\left(|c-\hat{c}|\right)\right]-\lambda\mathrm{CVaR}_\alpha\left((1-c)X\right)+\eta\mathrm{CVaR}^{1-\beta}\left((1-c)X\right), \tag{24}$$

where $g:[0,1]\to[0,+\infty)$ is twice the continuously differentiable, increasing, strictly convex, and Lipschitz continuous function with the Lipschitz constant $L_g>0$; the constant $\hat{c}\in[0,1]$ is the target diversion ratio; and the parameters $\lambda,\eta\geq 0$ are Lagrangian multipliers serving as weighting factors. In **Section 4**, we use $g\left(|c-\hat{c}|\right)=\frac{1}{2}(c-\hat{c})^2$.

The first term of (24) represents the diversion effort, where the deviation between the realized and targeted efforts is represented by the difference $|c-\hat{c}|$ and is measured through $g$. The second and third terms of (24) penalize too much and too less diversion of the flow, respectively. The target value $\hat{c}$ should be close to 1 for the hydropower generation case, as the water should be taken as much as possible, whereas it should be close to 0 for the flood diversion case because diverting too much water to the auxiliary channel should be avoided.

The second term of (24) is relevant for hydropower generation because extracting too much river water critically affects river aquatic ecosystems and leads to water rights conflicts among other industries, such as agriculture and manufacturing (Huuki et al., 2022; Basso et al., 2020; Widén et al., 2022). In contrast, the third term of (24) is relevant in the flood mitigation case, where high-flow events in a main channel should be diverted into an auxiliary channel to reduce the flood risk of downstream areas (Gopalan et al., 2022; Serra-Llobet et al., 2022). We discuss these two cases separately later. The theoretical analysis in this paper is based on the unified formulation (24), as it describes these cases simultaneously.

Set the optimization problem without model ambiguity:

$$V=\inf_{c\in\mathfrak{C}}J(c)=\inf_{c\in\mathfrak{C}}\left\{\mathbb{E}\left[g\left(|c(X)-\hat{c}|\right)\right]-\lambda\mathrm{CVaR}_\alpha\left((1-c(X))X\right)+\eta\mathrm{CVaR}^{1-\beta}\left((1-c(X))X\right)\right\}, \tag{25}$$

which is rewritten by (19)-(20) as follows:

$$V = \inf_{u,w\geq 0, c\in\mathfrak{C}} \left\{ -\lambda u + \eta w + \int \begin{pmatrix} g\left(\left|c(x)-\hat{c}\right|\right) + \dfrac{\lambda}{\alpha}\max\left\{u-\left(1-c(x)\right)x, 0\right\} \\ +\dfrac{\eta}{1-\beta}\max\left\{\left(1-c(x)\right)x - w, 0\right\} \end{pmatrix} p(x)\mathrm{d}x \right\}. \tag{26}$$

This is a convex optimization problem because the objective to be minimized is convex with respect to $u, w \geq 0, c \in \mathfrak{C}$, and their admissible set $[0,+\infty)^2 \times \mathfrak{C}$ is also convex. Therefore, $V$ exists.

### 3.2 Distributionally robust case

The distributionally robust optimization problem is formulated on the basis of the assumption that the PDF $p$ is not necessarily correct but that the true one (or the most reasonable one) is its transformation through a Radon–Nikodým derivative (Liu et al., 2022). Consider a family of PDFs $q_\omega$ parameterized by a positive measurable function $\omega = \omega(x)$ satisfying

$$q_\omega = \omega p \quad \text{in} \quad \Xi \tag{27}$$

and

$$\int_\Xi q_\omega(x)\mathrm{d}x = \int_\Xi \omega(x) p(x)\mathrm{d}x = 1. \tag{28}$$

The condition (27) means that the PDFs $p$ and $q_\omega$ are equivalent. The condition (28) is a natural requirement to ensure that $q_\omega$ is a PDF. Given $p$, the admissible set, which is the set of positive and measurable mappings $\omega$ satisfying the conditions (27)-(28), is denoted as $\Omega$. Choosing the constant function $\omega \equiv 1$ recovers $p = q_\omega$. Hereafter, we only consider $\omega$ belonging to $\Omega$. The integration range of integrands containing the PDFs $p, q_\omega$ are $\Xi$ and will be omitted in the integral representation for sake of simplicity.

Model ambiguity between the two PDFs $p$ and $q_\omega$ is evaluated by a divergence $D\left(p|q_\omega\right)$ as certain convex functional of $\omega$. Several divergences exist (Lin et al., 2022; Liu et al., 2022). We employ the relative entropy (Kullback–Leibler divergence) owing to its tractability:

$$D\left(p|q_\omega\right) = \int \phi\left(\omega(x)\right) p(x)\mathrm{d}x \tag{29}$$

with

$$\phi(x) = x\ln x - x + 1 \quad \text{for} \quad x > 0 \quad \text{and} \quad \phi(0) = 1. \tag{30}$$

Here, $\phi$ is non-negative and convex for $x \geq 0$ and is globally minimized at $x = 1$. Therefore, $D\left(p|q_\omega\right) \geq 0$ and it vanishes if and only if the two PDFs $p$ and $q_\omega$ are equivalent. The two PDFs $p$ and $q_\omega$ are more different from each other for larger $D\left(p|q_\omega\right)$. We demonstrate that Kullback–Leibler divergence is able to mathematically handle the unbounded state (streamflow discharge) under uncertainties.

The distributionally robust formulation of the stochastic optimization problem is

$$W=\inf_{u,w\geq 0,c\in\mathfrak{C}}\sup_{\omega\in\Omega}\left\{\begin{array}{l}-\lambda u+\eta w\\+\int\left(\begin{array}{l}g\left(\left|c(x)-\hat{c}\right|\right)+\dfrac{\lambda}{\alpha}\max\left\{u-\left(1-c(x)\right)x,0\right\}\\+\dfrac{\eta}{1-\beta}\max\left\{\left(1-c(x)\right)x-w,0\right\}\end{array}\right)q_{\omega}(x)\mathrm{d}x-\mu D\left(p|q_{\omega}\right)\end{array}\right\} \quad (31)$$

with a Lagrangian multiplier $\mu>0$. Now, the expectations are evaluated under the distorted PDF $q_{\omega}$ but not with $p$. The last term (31) works as a penalization of the model ambiguity. The decision maker becomes more ambiguity averse for smaller $\mu$.

Hereafter, we always employ **Assumption 1** below.

***Assumption 1*** $C_X\equiv\dfrac{2\eta}{\mu(1-\beta)}>0$ *is sufficiently small such that* $\mathbb{E}\left[X^k e^{C_X X}\right]<+\infty$ ($k=0,1,2$).

This is possible both for supOU and supCBI cases because of **Propositions 1-3** if $B$ is sufficiently small.

***Remark 3*** One may evaluate the diversion effort by $g\left(\left|c-\hat{c}\right|x\right)$, which evaluates the deviation between the realized and targeted discharges but not the diversion ratios. This modification does not fit into our theory because one may additionally need $\mathbb{E}\left[e^{X^{\sigma}}\right]<+\infty$ ($\sigma>0$).

The right-hand side of (31) can be rewritten in a more concise form that is not only useful in applications but also provides better theoretical insights into the distributionally robust optimization problem. The inner optimization problem for $\omega$ is explicitly solved as follows. Given $u,w\geq 0,c\in\mathfrak{C}$, write

$$F(x,u,w,c)=g\left(\left|c(x)-\hat{c}\right|\right)+\frac{\lambda}{\alpha}\max\left\{u-\left(1-c(x)\right)x,0\right\}+\frac{\eta}{1-\beta}\max\left\{\left(1-c(x)\right)x-w,0\right\}. \quad (32)$$

This $F$ is nonnegative. We obtain

$$\begin{aligned}&\int\left(\begin{array}{l}g\left(\left|c(X)-\hat{c}\right|\right)+\dfrac{\lambda}{\alpha}\max\left\{u-\left(1-c(x)\right)x,0\right\}\\+\dfrac{\eta}{1-\beta}\max\left\{\left(1-c(x)\right)x-w,0\right\}\end{array}\right)q_{\omega}(x)\mathrm{d}x-\mu D\left(p|q_{\omega}\right)\\&=\int\left\{F(x,u,w,c)\,\omega(x)-\mu\phi\left(\omega(x)\right)\right\}p(x)\mathrm{d}x\end{aligned}. \quad (33)$$

Since $p\geq 0$, the quantity on the right-hand side subject to $\omega\in\Omega$ is minimized by a constrained optimization (e.g., Hansen and Miao (2018)):

$$\omega^{*}(x)=\frac{e^{\theta(x,u,w,c)}}{\int e^{\theta(x,u,w,c)}p(x)\mathrm{d}x}\equiv\frac{1}{K}e^{\theta(x,u,w,c)}\ \text{ with }\ \theta(x,u,w,c)=\frac{F(x,u,w,c)}{\mu}. \quad (34)$$

Note that $\omega^{*}\in\Omega$ for each $u,w\geq 0,c\in\mathfrak{C}$. Therefore, the inner maximization of (31) is solved as

follows:

$$\begin{aligned}
&\sup_{\omega\in\Omega}\left\{-\lambda u+\eta w+\int F(x,u,w,c)q_{\omega}(x)\mathrm{d}x-\mu D(p|q_{\omega})\right\} \\
&=-\lambda u+\eta w+\int\left\{F(x,u,w,c)\omega^{*}(x)-\mu\left(\omega^{*}(x)\right)\right\}p(x)\mathrm{d}x\,. \\
&=-\lambda u+\eta w+\mu\int\left(\ln K+\frac{1}{K}e^{\theta(x,u,w,c)}-1\right)p(x)\mathrm{d}x
\end{aligned} \tag{35}$$

Furthermore, we have

$$\int\left(\ln K+\frac{1}{K}e^{\theta(X,u,w,c)}-1\right)p(x)\mathrm{d}x=\ln K+\frac{1}{K}K-1=\ln\left(\int e^{\frac{F(x,u,w,c)}{\mu}}p(x)\mathrm{d}x\right). \tag{36}$$

Consequently, we arrive at

$$W=\inf_{u,w\geq 0,c\in\mathfrak{C}}\left\{-\lambda u+\eta w+\mu\ln\left(\int e^{\frac{F(x,u,w,c)}{\mu}}p(x)\mathrm{d}x\right)\right\}\equiv\inf_{u,w\geq 0,c\in\mathfrak{C}}\left\{-\lambda u+\eta w+H(u,w,c)\right\} \tag{37}$$

with $H(u,w,c)=\mu\ln\left(\int e^{\frac{F(x,u,w,c)}{\mu}}p(x)\mathrm{d}x\right)$, showing that (31) is an optimization problem where the certainty equivalent $\mu\ln\mathbb{E}\left[\exp\left(\frac{(\cdot)}{\mu}\right)\right]$ of $\mathbb{E}$ is employed. We have $V\leq W$ by the classical Jensen's inequality. In addition, we have $W|_{\mu=\mu_1}\leq W|_{\mu=\mu_2}$ for $0<\mu_2\leq\mu_1$ (**Lemma A.1** of **Appendix A**), showing that the smaller penalization against the ambiguity leads to a larger (i.e., worse) estimate of the objective.

The admissible range $[0,+\infty)^2$ of the two decision variables $u,w$ is unbounded. However, we can limit our interest to strictly bounded ones, as shown in **Proposition 4** below. In this way, the admissible set of decision variables $u,w$ can be assumed to be compact *a priori*.

***Proposition 4*** *Any minimizing $u,w$ of (37) belong to some compact set $D\subset[0,+\infty)^2$.*

***Remark 4*** Essentially, the same compactification argument applies to the regularized optimization problem and its discretized counterpart in the next subsection. We can choose a sufficiently large compact set $D$ to which the optimal decision variables of the three optimization problems (the distributionally robust one and its regularization with/without discretization in the next subsection) belong. This compaction is introduced for purely theoretical reasons.

By **Proposition 4**, we will often use a suitably wide and compact set $D$ for the decision variables $u,w$ without loss of generality. We have

$$
\begin{aligned}
F(x,u,w,c) &\le g(1)+\frac{\lambda}{\alpha}\max\{u-(1-c(x))x,0\}+\frac{\eta}{1-\beta}\max\{(1-c(x))x-w,0\} \\
&\le g(1)+\frac{\lambda u}{\alpha}+\frac{\eta}{1-\beta}(1-c(x))x \\
&\le g(1)+\frac{\lambda U}{\alpha}+\frac{\eta}{1-\beta}x \\
&\equiv \bar{F}+\frac{\eta}{1-\beta}x
\end{aligned}
\quad , \ (u,w)\in D, c\in\mathfrak{C} \qquad (38)
$$

with a constant $U>0$, which exists by **Proposition 4**.

To deeper investigate the optimization problem, we use the lower semi-continuity result.

***Proposition5*** *$H=H(u,w,c)$ is lower semi-continuous in $[0,+\infty)^2\times\mathfrak{C}$.*

Because $-\lambda u+\eta w+H_\tau(u,w,c)$ is lower semi-continuous, the optimization problem (37) admits a solution if it is convex. Set

$$G(\theta)=\ln\left(\int e^{\theta(x)}p(x)\mathrm{d}x\right) \qquad (39)$$

for any function $\theta:\Xi\to[0,+\infty)$ such that $\int e^{\theta(x)}p(x)\mathrm{d}x<+\infty$.

***Proposition 6*** *For $\theta=\theta_1,\theta_2$ with $\int e^{\theta_i(x)}p(x)\mathrm{d}x<+\infty$ ($i=1,2$) and any real number $\kappa\in[0,1]$*

$$G(\kappa\theta_1+(1-\kappa)\theta_2)\le\kappa G(\theta_1)+(1-\kappa)G(\theta_2). \qquad (40)$$

*Furthermore, for $u_i,w_i\ge 0, c_i\in\mathfrak{C}$ ($i=1,2$),*

$$H(\kappa u_1+(1-\kappa)u_2,\kappa w_1+(1-\kappa)w_2,\kappa c_1+(1-\kappa)c_2)\le\kappa H(u_1,w_1,c_1)+(1-\kappa)H(u_2,w_2,c_2). \qquad (41)$$

*Consequently, the optimization problem (37) is convex.*

The convexity (**Proposition 6)** and the lower semi-continuity (**Proposition 5**) of $H$ imply that it is sequentially weakly lower semi-continuous (Theorem 2.20 of Galewski (2021)) and that the problem has a globally minimizing $(u,w,c)$ with which the infimum is attained.

### 3.3 Regularization

Our distributionally robust optimization problem is non-smooth term as it has "max" functions in the objective. This non-smoothness prevents us from applying a gradient descent type solution algorithm to the optimization problem. This issue is resolved through the replacement of the "max" function by a smoother one (e.g., Luna et al. (2016)). We use the strictly convex regularized "max" function as in Yoshioka et al. (2023c):

$$m_\tau(x)=\frac{x+\sqrt{x^2+4\tau^2}}{2} \quad \text{for} \quad x\in\mathbb{R} \tag{42}$$

with a regularization parameter $\tau>0$. We note a key estimate used in the sequel:

$$\max\{x,0\}\le m_\tau(x)\le \max\{x,0\}+\tau \quad \text{for} \quad x\in\mathbb{R}, \tag{43}$$

showing that $m_\tau(\cdot)$ uniformly approximates $\max\{\cdot,0\}$.

For $\tau>0$, the regularized distributionally robust optimization problem is formulated as follows:

$$W_\tau=\inf_{u,w\ge 0,c\in\mathfrak{C}}\left\{-\lambda u+\eta w+\mu\ln\left(\int e^{\frac{F_\tau(x,u,w,c)}{\mu}}p(x)\,\mathrm{d}x\right)\right\}\equiv\inf_{u,w\ge 0,c\in\mathfrak{C}}\left\{-\lambda u+\eta w+H_\tau(u,w,c)\right\}, \tag{44}$$

where

$$F_\tau(x,u,w,c)=g(|c-\hat{c}|)+\frac{\lambda}{\alpha}m_\tau(u-(1-c)x)+\frac{\eta}{1-\beta}m_\tau((1-c)x-w). \tag{45}$$

We have that $-\lambda u+\eta w+H_\tau(u,w,c)$ is strictly convex.

***Proposition 7*** *For any* $\kappa\in(0,1)$, $u_i,w_i\ge 0,c_i\in\mathfrak{C}$ *(* $i=1,2$ *) with* $(u_1,w_1,c_1)\ne(u_2,w_2,c_2)$,

$$H_\tau(\kappa u_1+(1-\kappa)u_2,\kappa w_1+(1-\kappa)w_2,\kappa c_1+(1-\kappa)c_2)<\kappa H_\tau(u_1,w_1,c_1)+(1-\kappa)H_\tau(u_2,w_2,c_2). \tag{46}$$

*Thereby, the optimization problem (44) has a strictly convex objective and admits a unique minimizer in* $D\times\mathfrak{C}$ *where* $D\subset[0,+\infty)^2$ *is a compact set.*

The optimal, i.e., worst-case ambiguity $\omega=\omega^*$ is found as

$$\omega^*(x)=\frac{e^{\frac{F_\tau(x,u^*,w^*,c^*(x))}{\mu}}}{\int e^{\frac{F_\tau(x,u^*,w^*,c^*(x))}{\mu}}p(x)\,\mathrm{d}x} \tag{47}$$

with $(u^*,w^*,c^*)$ the unique minimizer of (44).

Finally, for consistency between the regularized and nonregularized optimization problems, we propose the following.

***Proposition 8*** *It follows that* $W_\tau\to W$ *as* $\tau\to+0$. *In addition, the minimizer of* $W_\tau$ *accumulates to a minimizer of* $W$. *Furthermore,* $|W_\tau-W|\le C\tau$ *with a constant* $C>0$ *independent from* $\tau$.

### 3.4 Discretization of the regularized problem

The discretized counterpart of (44) is defined as follows:

$$W_{N,\tau}=\inf_{u,w\ge 0,c\in\mathfrak{C}_N}\left\{-\lambda u+\eta w+\mu\ln\left(\sum_{i=1}^N p_i e^{\frac{F_{i,\tau}(u,w,c)}{\mu}}\right)\right\} \tag{48}$$

with a discrete probability measure $P_N(\mathrm{d}X)=\sum_{i=1}^{N} p_i \delta_{(X=X_i)}$ ($\sum_{i=1}^{N} p_i = 1$, $p_i > 0$, $0 < X_1 < X_2 < ... < X_n$ for $1 \le i \le N$) approximating $p$, the coefficients

$$F_{i,\tau}(u,w,c) = g(|c_i - \hat{c}|) + \frac{\lambda}{\alpha}\max\{u - (1-c_i)X_i, 0\} + \frac{\eta}{1-\beta}\max\{(1-c_i)X_i - w, 0\}, \tag{49}$$

and the admissible set of decision variables

$$\mathfrak{C}_N = \left\{ c = \{c_i\}_{1\le i\le N} \,\middle|\, c_i \in [0,1],\ 1 \le i \le N \right\}. \tag{50}$$

The discretized worst-case ambiguity $\omega_N$ is

$$\omega^* = \{\omega_i^*\}_{1\le i\le N} \quad \text{with} \quad \omega_i^* = \frac{e^{\frac{F_{i,\tau}(u,w,c)}{\mu}}}{\sum_{j=1}^{N} p_j e^{\frac{F_{j,\tau}(u,w,c)}{\mu}}}. \tag{51}$$

The strict convexity and hence the unique solvability of the discretized version of the regularized problem follow in essentially the same way as the continuous one, which is therefore omitted.

We end this section by providing comments on the convergence from discretized to continuous problems. We expect that a reasonable discretization should satisfy $W_{N,\tau} \to W_\tau$ under $N \to +\infty$ for each fixed $\tau > 0$. Usually, the discretization of probability densities relies on the continuous mapping theorem (Theorem 2.12 of Jiang (2022)):

$$\lim_{N\to+\infty} \int f(x) P_N(\mathrm{d}x) = \int f(x) P(\mathrm{d}x) \quad \text{for any} \quad f \in C(\mathbb{R}), \tag{52}$$

the assumption below is assumed to be satisfied by $P_N$.

***Assumption 2*** $P_N(\mathrm{d}x)$ *weakly converges to* $P(\mathrm{d}x) = p(x)\mathrm{d}x$.

This assumption is satisfied for the discretization method of the probability density (Yoshioka et al., 2023a) employed in our application. To ensure that **Assumption 2** is satisfied, we need that $F_\tau(x,u,w,c(x))$ is continuous for $x > 0$ under optimality, which is difficult to judge *a priori*. Nevertheless, the continuity assumption of $c$ will be experimentally checked in our numerical computation, as shown in the next section. If $F_\tau(x,u,w,c(x))$ is continuous for $x > 0$, then under a suitable Hilbert space, we will have convergence, up to the subsequence, $W_{N,\tau} \to W_\tau$ as $N \to +\infty$ if $J(u,w,c) \le J_N(u,w,c) + O_N$ ($(u,w,c) \in D \times \mathfrak{C}$, $N \in \mathbb{N}$) with some $O_N$ independent of $(u,w,c) \in D \times \mathfrak{C}$ such that $\liminf_{N\to+\infty} O_N = 0$.

### 3.5 Numerical algorithm

We present a numerical algorithm to solve the discretized version of the regularized distributionally robust optimization problem (48). Because this paper is not aiming at comparing optimization methods but rather addressing a case study, it suffices to develop a simple algorithm for its purpose. It should be noted that this point represents a limitation of this paper as a case study. Our algorithm is a gradient descent method with an inertia, which has successfully been applied to a variety of problems (Calder and Yezzi, 2019; Duruisseaux and Leok, 2022; He et al., 2022; Wang et al., 2022). Set

$$\hat{H}\left(u,w,c\right)=-\lambda u+\eta w+\mu \ln\left(\sum_{i=1}^{N} p_i e^{\frac{F_{i,\tau}(u,w,c)}{\mu}}\right). \tag{53}$$

In **Algorithm 1** below, the superscript represents the iteration count. Parameters $\Delta t\left(=O\left(10^{-4}\right)\right)$, $\gamma\left(=O\left(10^{-1}\right)\right)$, and $\delta=O\left(10^{0}\right)$ control the convergence speed. The parameter $\varepsilon\left(=10^{-8}\right)$ is an error tolerance.

***Algorithm 1***

***Step 0*** *Prepare model parameter values and set the initial guesses* $\left(u^{(0)},w^{(0)},c^{(0)}\right)=\left(0,0,\{0\}_{1\le i\le N}\right)$ *and* $\left(u^{(-1)},w^{(-1)},c^{(-1)}\right)=\left(u^{(0)},w^{(0)},c^{(0)}\right)$. *Set* $n=0$.

***Step 1*** *Compute* $\left(u^{(n+1)},w^{(n+1)},c^{(n+1)}\right)$ *as follows:*

$$u^{(n+1)}=\max\left\{u^{(n)}-\frac{\delta\Delta t}{\lambda}\frac{\partial\hat{H}(u,w,c)}{\partial u}\bigg|_{(u,w,c)=\left(u^{(n)},w^{(n)},c^{(n)}\right)}+\gamma\left(u^{(n)}-u^{(n-1)}\right),0\right\}, \tag{54}$$

$$w^{(n+1)}=\max\left\{v^{(n)}-\frac{\delta\Delta t}{\eta}\frac{\partial\hat{H}(u,w,c)}{\partial w}\bigg|_{(u,w,c)=\left(u^{(n)},w^{(n)},c^{(n)}\right)}+\gamma\left(w^{(n)}-w^{(n-1)}\right),0\right\}, \tag{55}$$

$$c_i^{(n+1)}=\min\left\{\max\left\{c_i^{(n)}-\frac{\Delta t}{p_i}\frac{\partial\hat{H}(u,w,c)}{\partial c_i}\bigg|_{(u,w,c)=\left(u^{(n)},w^{(n)},c^{(n)}\right)}+\gamma\left(c_i^{(n)}-c_i^{(n-1)}\right),0\right\},1\right\}\quad (1\le i\le N). \tag{56}$$

***Step 2*** *Compute the update error as follows:*

$$e^{(n)}=\min\left\{\left|u^{(n+1)}-u^{(n)}\right|,\left|w^{(n+1)}-w^{(n)}\right|,\min_{1\le i\le N}\left\{\left|c_i^{(n+1)}-c_i^{(n)}\right|\right\}\right\}\quad (1\le i\le N). \tag{57}$$

*If* $e^{(n)}\le\varepsilon$, *then output* $\left(u^{(n+1)},w^{(n+1)},c^{(n+1)}\right)$ *as the optimal decision variables* $\left(u^{*},w^{*},c^{*}\right)$ *and compute the worst-case ambiguity by (51) based on* $\left(u^{*},w^{*},c^{*}\right)$. *If* $e^{(n)}>\varepsilon$, *then update* $n\to n+1$ *and go to* ***Step 1***.

**Algorithm 1** is a gradient descent method with an inertia, where the inertial effect to accelerate the convergence of the optimization problem is implemented with the terms multiplied by $\gamma$. The strict convexity of the objective ensures the convergence of the algorithm to the unique minimizer if $\Delta t$ is small and $\gamma$ and $\tau$ are chosen properly. Some researchers have analyzed time-dependent (iteration count-dependent) $\gamma$ (Ushiyama et al., 2022), but we consider a simpler one with a constant $\gamma$. This simplest setting significantly improves the pure gradient descent ($\gamma = 0$) (**Figures 11-12** shown later). Projections using "min" and "max" functions are employed in (54), (55), and (56) to guarantee that the computed decision variables rigorously belong to their admissible sets. This projection is understood as a large-penalty limit of the implicit penalization method (Calder and Yezzi, 2019). The parameter $\delta$ appearing in (54)-(55) is a scaling parameter to accelerate the convergence that has been chosen manually. The division by $p_i$ in (56) is to accelerate the computation for small probabilities.

## 4. Applications

The study sites are in the Hii River system in the San-in area of Shimane Prefecture, Japan, where we have been investigating the local hydrology and ecology since 2014. **Figure 2** shows the two study sites, which have different hydrological characteristics and face different optimization issues. Throughout this section "$X$" stands for the discharge value but not a stochastic process. We formulate the problem within the context of our theoretical analysis for modeling simplicity, although there may include some nonconvexity in the other engineering applications. Here, **Assumption 1** is satisfied in our MMA processes. **Assumption 2** is a technical mathematical assumption. The strict convexity is also a technical assumption to formulate the optimization problem, which is the limitation of this paper.

### 4.1 Model identification and Markovian lift

The PDF $p$ of the target process is computed as follows. We do not present a detailed explanation of the Markovian lift and Fourier transform in this paper, but such information can be found in the literature (Yoshioka et al., 2023a; Hainaut, 2022). The degrees of freedom $n$ of the Markovian lift is fixed to 1,024 according to the quantile discretization (Yoshioka et al., 2023a), which has preliminary been found to be sufficient for the computation below. The degrees of freedom $N$ for the discretization of the optimization problem are set to 2,000 for **Case 1** (i.e., control of the low flow) and 8,000 for **Case 2** (i.e., control of the high flow). The domains of the cases have lengths of 200 m/s and 1,000 $m^3/s$, respectively. **Appendix B** summarizes the statistics of the discharge time series at each study site and presents the identified model parameter values of the supCBI process.

We set the regularization parameter $\tau = 0.0001$ $m^3/s$ for **Case 1** and $\tau = 0.01$ $m^3/s$ for **Case 2**, which have been found to be sufficiently small as well as computationally stable. Larger $\tau$ would be necessary for problems with larger mean discharge and/or large increment $\Delta t$ in **Algorithm 1**.

### 4.2 Case 1: Control of the low flow

#### 4.2.1 Study site

Water allocation for hydropower generation is a crucial issue for the main channel of the Hii River upstream of Obara Dam, which is the largest dam in the river (**Figure 3**). There is a hydropower plant with an intake point upstream of the reservoir created by Obara Dam that has been idling from June 2020 to probably summer or autumn 2023 (Yoshioka and Yoshioka, 2022a). The hydropower plant is expected to provide a source of clean electricity for the regional grid, but an intake of too much water will affect the downstream dam–reservoir system. The study site is a habitat of the landlocked *Plecoglossus altivelis altivelis* (i.e., ayu sweetfish) and thus serves as a major inland fishery resource. The spawning and migration of this fish are negatively affected by extremely low-flow regimes (Maeda, 2013). In addition, the intake point may generate a high-flow reach that can potentially entrain fish resulting in death or injury by impingement at the screen (Jarvis and Closs, 2019). Therefore, balancing the water intake policy for future operation of the hydropower plant is a crucial issue. **Appendix B** shows that the streamflow discharge at this study site has a long memory, and hence the application of the MMA process is relevant.

**Figure 4** compares the observed and modeled PDFs of the streamflow discharge in ordinary and logarithmic scales. The tails of the PDFs agree well with each other by **Figure 4(b),** but they are visibly different near the modes because the modeled PDF is unimodal while the observed PDF is seemingly bimodal. However, the observed PDF at the study site is unimodal before May 2020 (Yoshioka and Yoshioka, 2022a). The bimodality of a discharge PDF is often due to time-periodic fluctuations (Ye et al., 2021). We attributed the bimodality in this case to the hydropower generation by an upstream smaller dam aimed at compensating for the lack of electricity supply from the idling plant. The deviation between the observation and model implies that solving the optimization problem under ambiguity is essential.

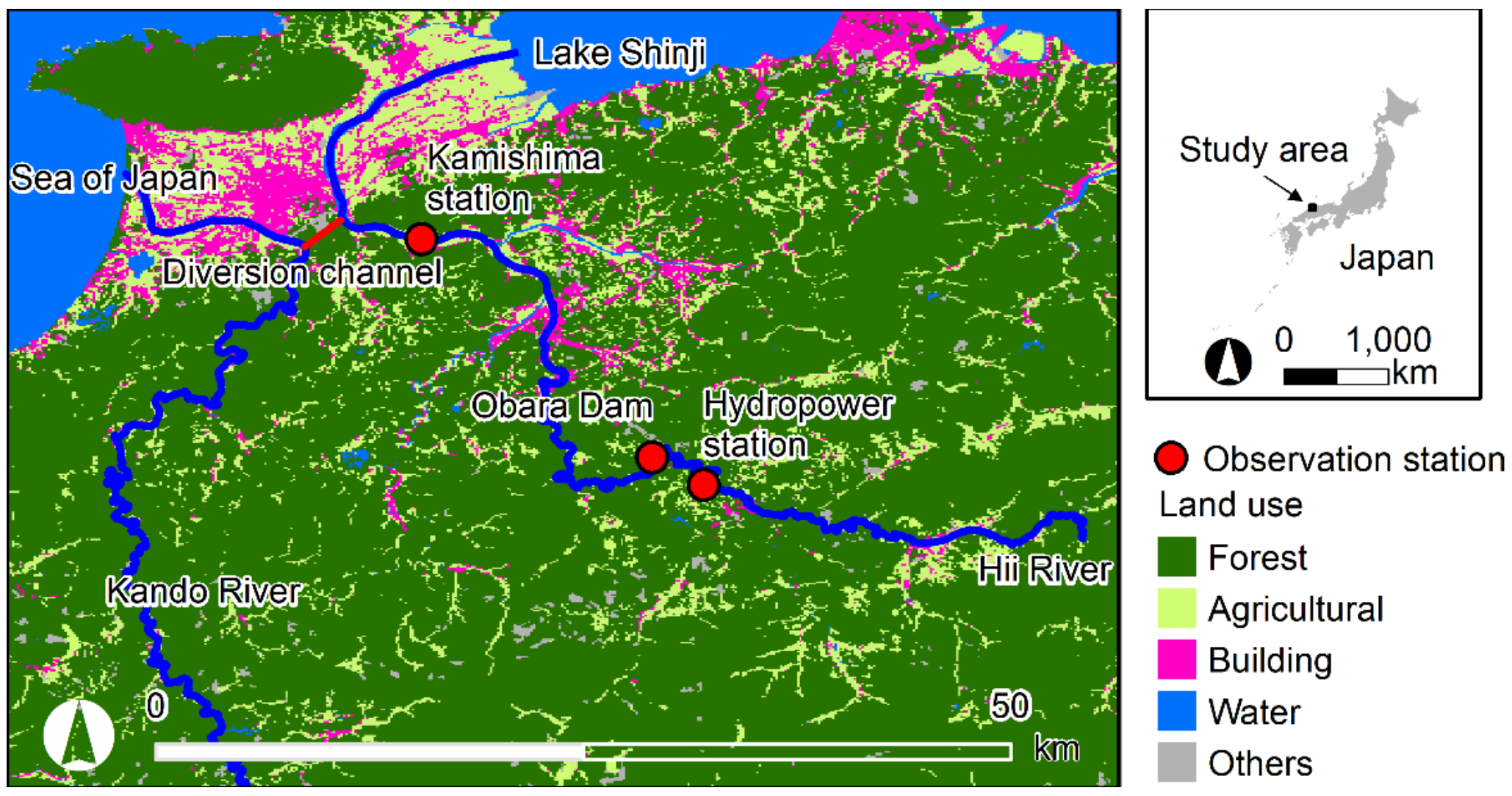

**Figure 2**. Map of the study area, including the hydropower station (study site of **Case 1**) and diversion channel (study site of **Case 2**).

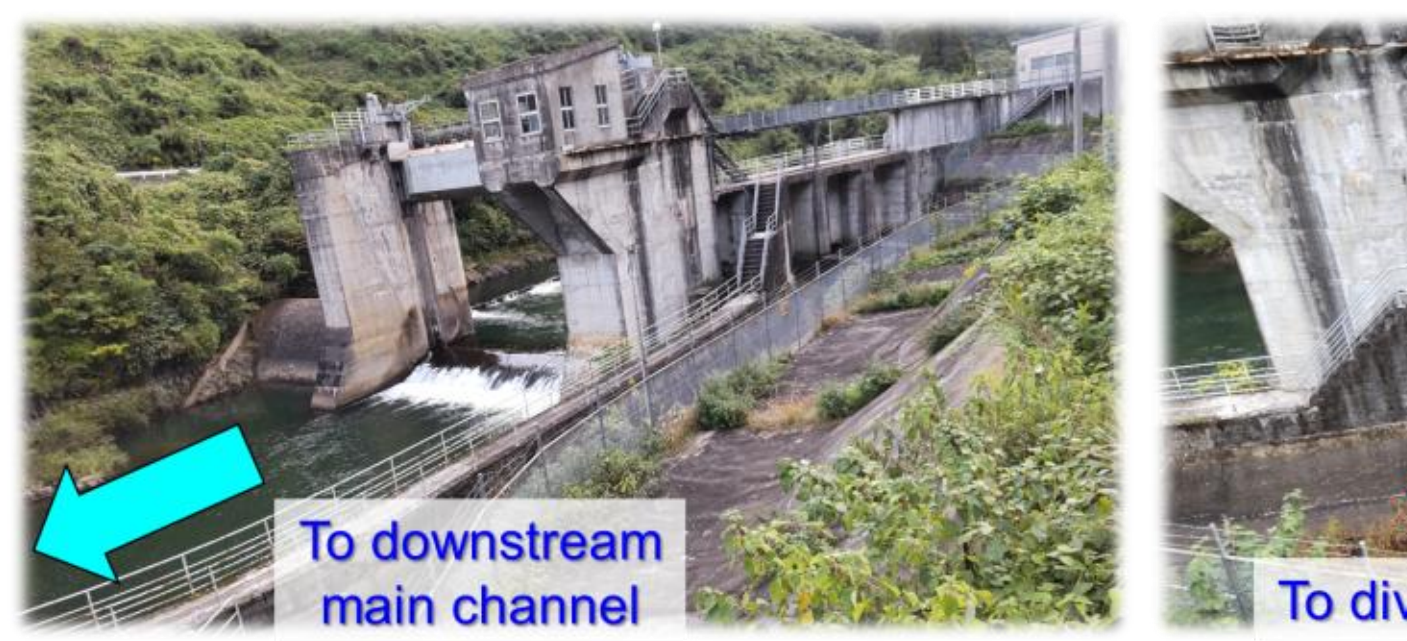

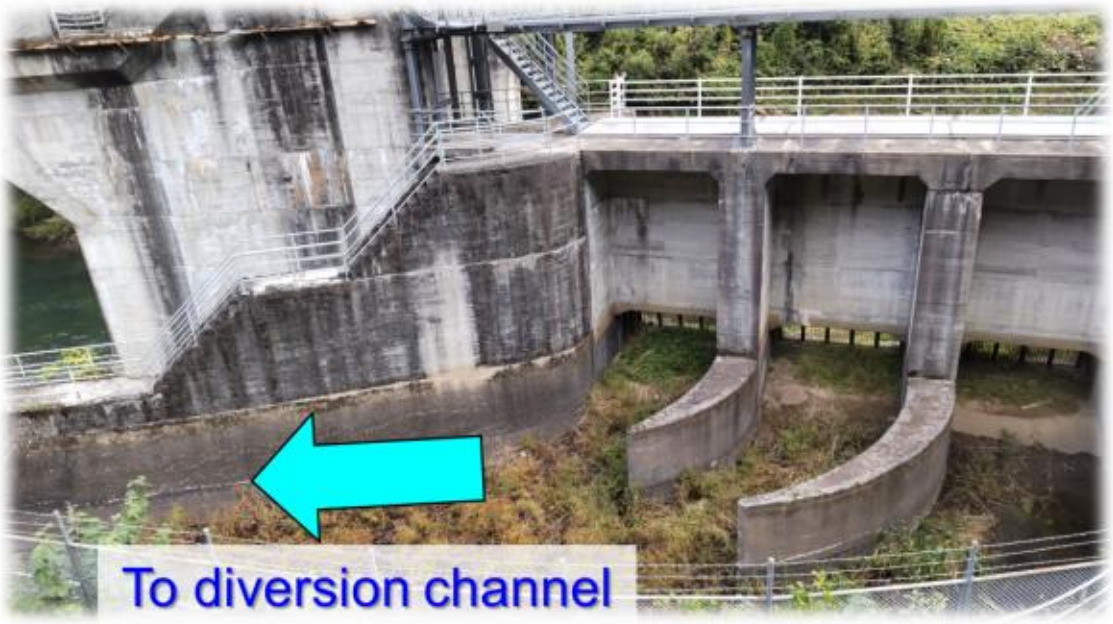

**Figure 3.** Photographs of the intaking point of the hydropower station (a) taken from the downstream and (b) from the diversion channel. Photo taken by Yoshioka H. on October 19, 2022.

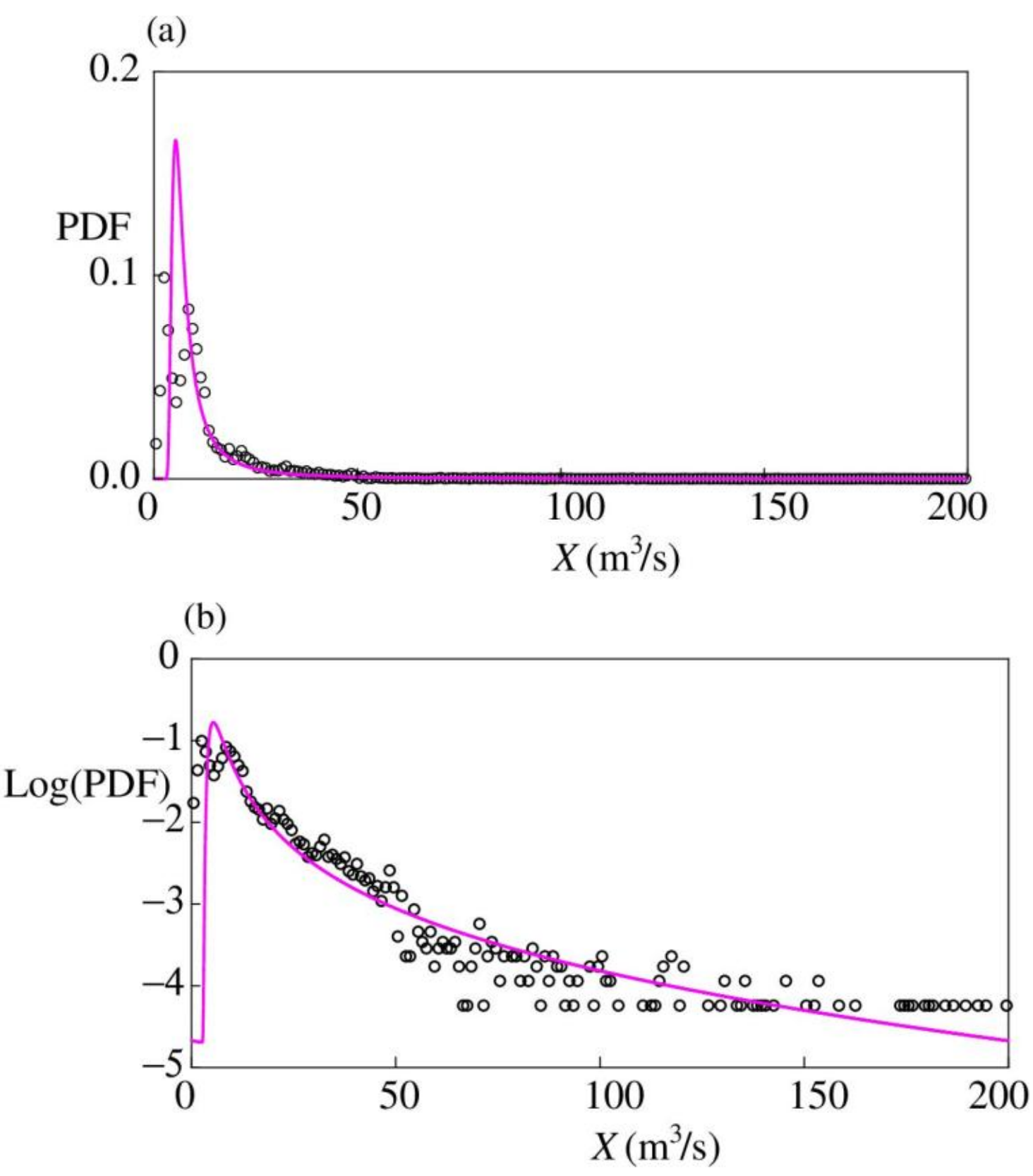

**Figure 4**. Observed and modelled PDFs on **(a)** ordinary scale and **(b)** logarithmic scale.

#### 4.2.2 Computation and sensitivity analysis

**Figures 5**, **6**, and **7** show the computed diversion ratio $c$, ambiguity $\omega$, and PDF $q$ under the optimality, respectively, as bivariate functions of the discharge $X$ and weight $\lambda$ under the worst-case optimization. The other parameters were $\alpha = 0.20$ and $\mu = 1.0$. The computational results are presented for a small discharge because the objective of **Case 1** is the optimization of low-flow events. **Figure 5** shows that the optimized diversion ratio depends nonmonotonically on both $X$ and $\lambda$, but the dependence seems to be smooth except for $X = 0$ where $c$ is computed to be very close to 1 numerically. However, this is not problematic because no water can be abstracted at a zero discharge $X = 0$. The river water should not be abstracted for a small discharge, which is consistent with the intuition that a persistent river flow should be kept for sustainable water resources management. Interestingly, the proposed model suggests that expanding the area of no-intake with increasing ambiguity aversion is suboptimal. **Figures 6** and **7** show that the impact of model ambiguity is limited to small discharges at both small and large $\lambda$, which can be attributed to the use of a lower CVaR for risk evaluation.

We also investigate the sensitivity against another weight $\mu$ for ambiguity aversion. **Figures 8**, **9**, and **10** show the computed diversion ratio $c$, ambiguity $\omega$, and PDF $q$ under the optimality, respectively, as bivariate functions of the discharge $X$ and weight $\mu$ under the worst-case optimization. The other parameters were $\alpha = 0.20$ and $\lambda = 1$. As shown in **Figure 8**, the proposed model finds that expanding the area of no-intake was only optimal with increasing ambiguity aversion (i.e., decreasing $\mu$ ). Even abstracting very low flows close to $X = 0$ are optimal for small $\mu$ such that the ambiguity and PDF were concentrated near $X = 0$. As shown in **Figures 9** and **10**, the dependence of the ambiguity and PDF on ambiguity aversion was nonmonotonic, which suggests that the benchmark PDF should be distorted gradually as $\mu$ decreases to $O(10^{-1})$. A risk-averse decision-maker should consider a PDF significantly concentrating on $X = 0$.

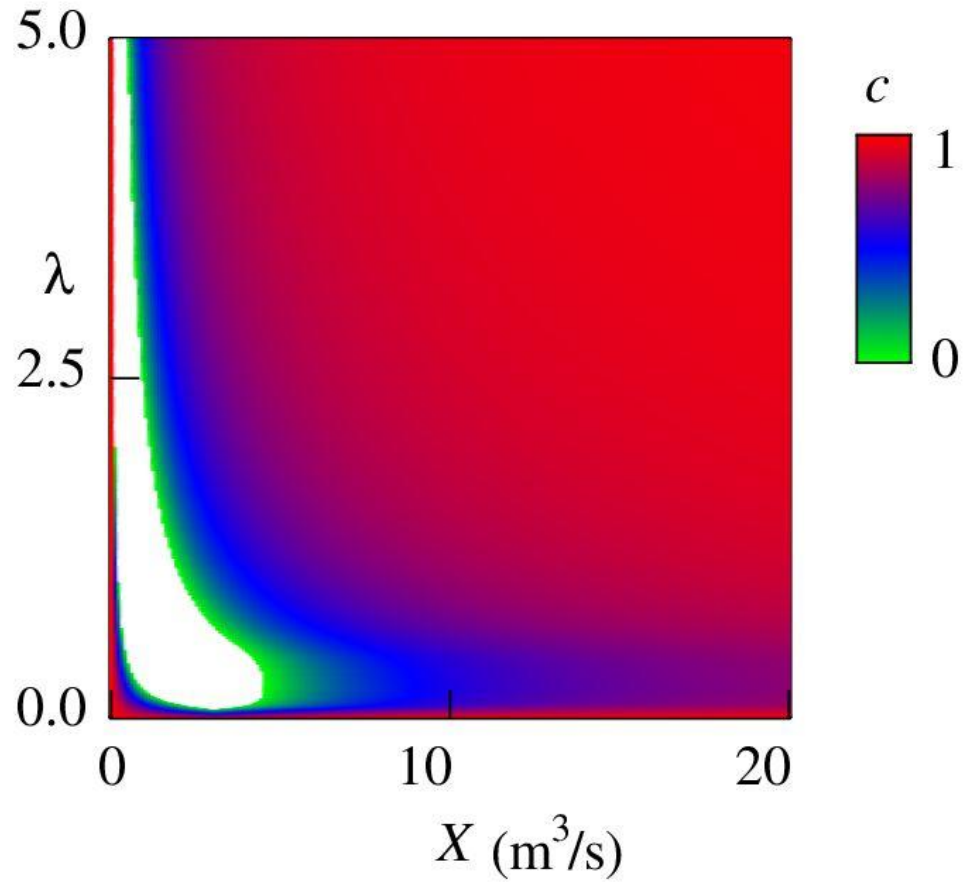

**Figure 5.** The diversion ratio under the worst-case optimization as a function of $(X, \lambda)$. The white color is the area of $c = 0$ in which the water should not be abstracted.

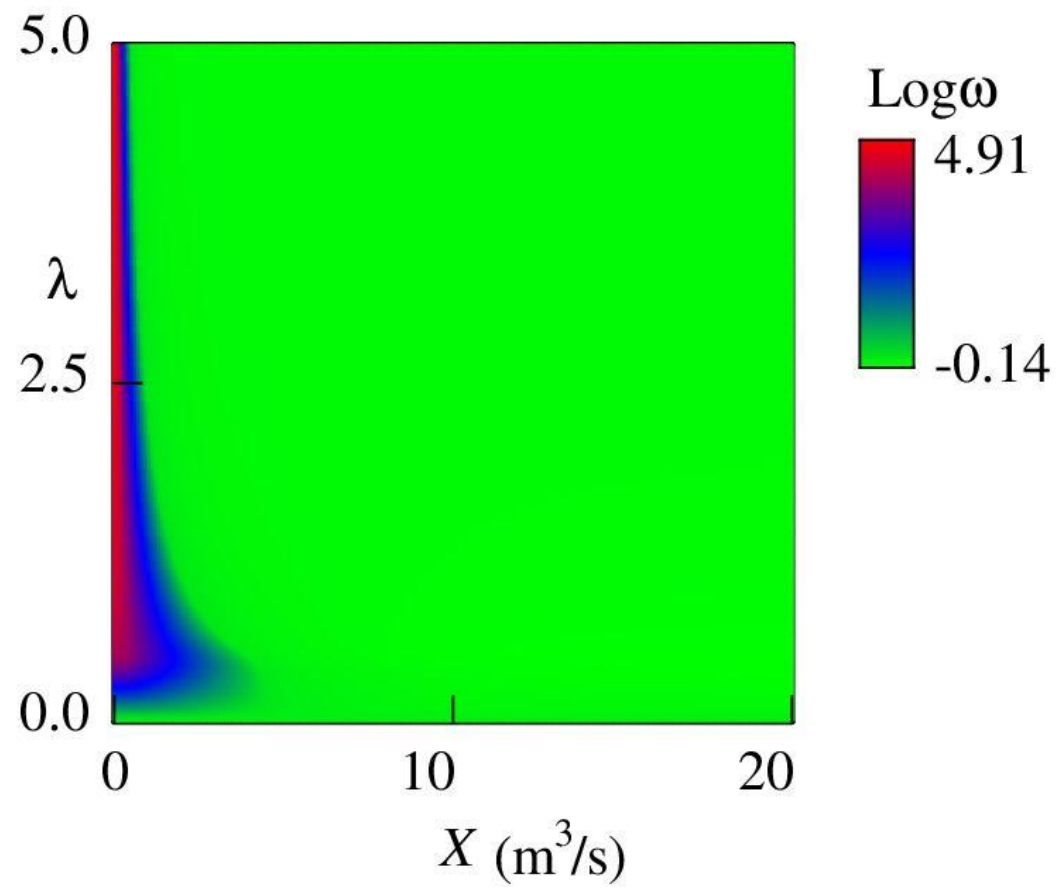

**Figure 6.** The ambiguity under the worst-case optimization as a function of $(X, \lambda)$.

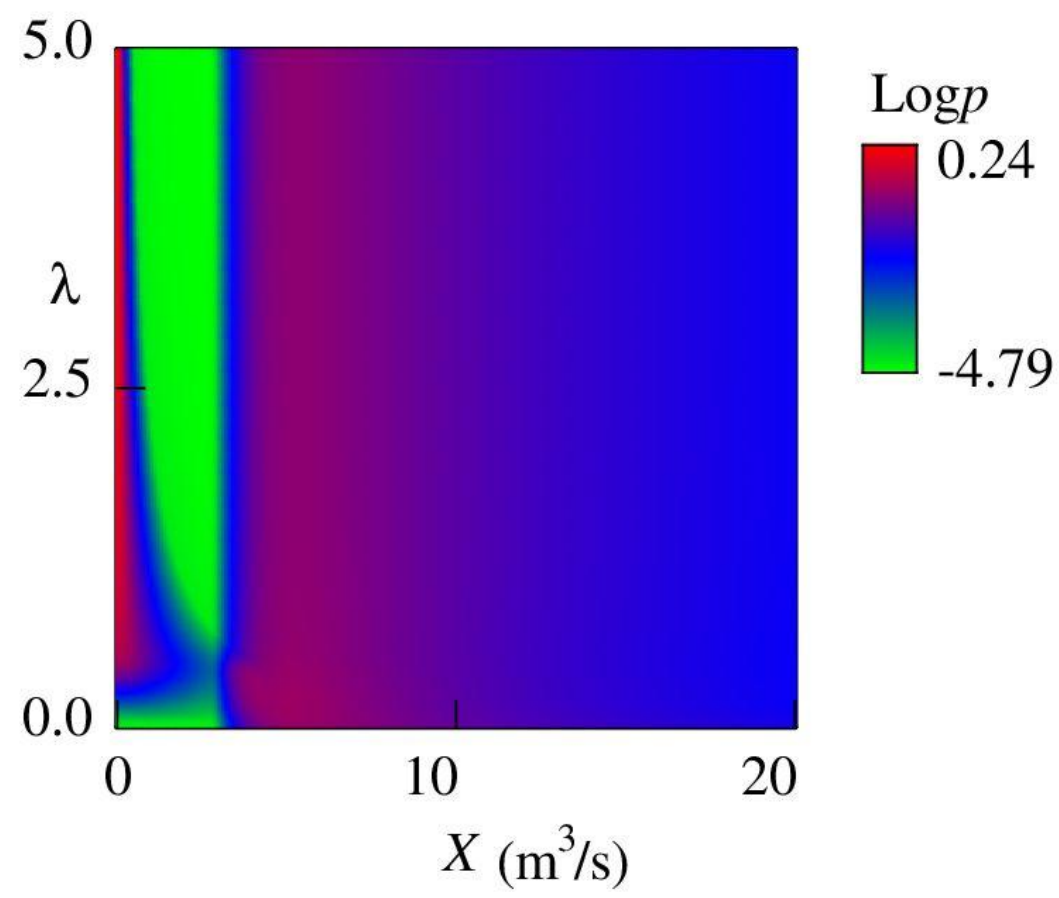

**Figure 7.** The PDF under the worst-case optimization as a function of $(X, \lambda)$.

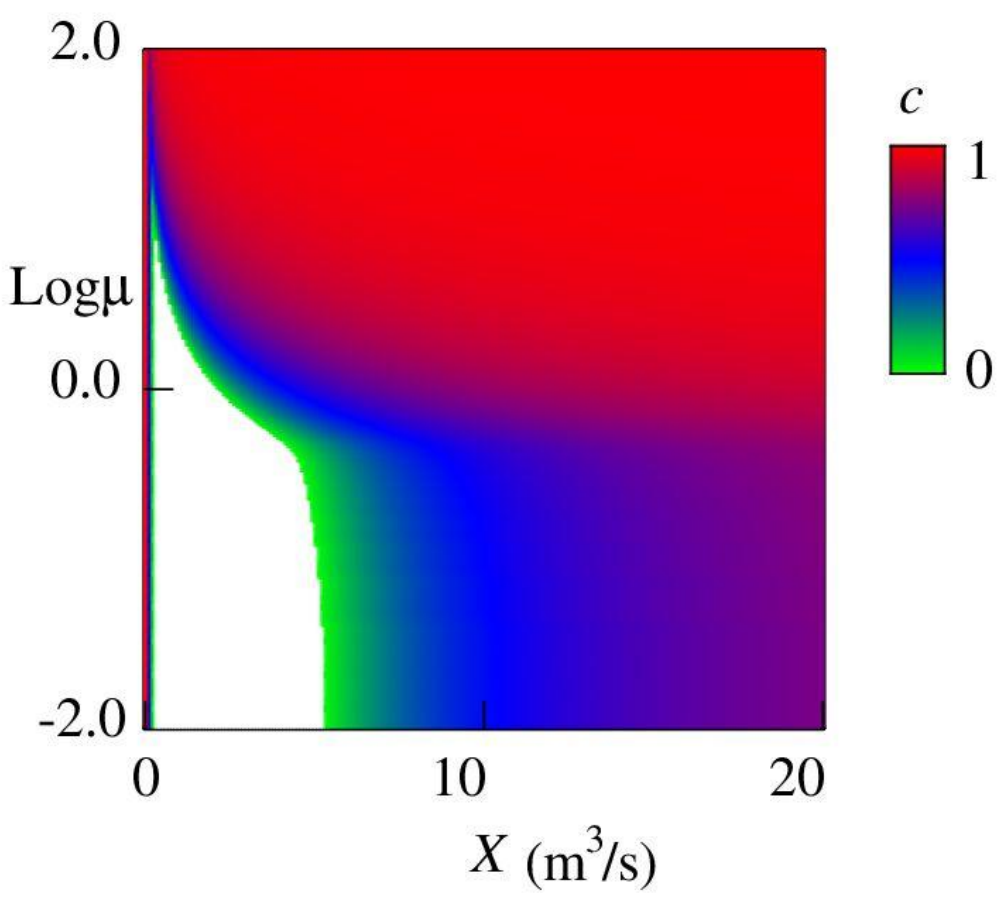

**Figure 8.** The diversion ratio under the worst-case optimization as a function of $(X, \mu)$. The white color is the area of $c = 0$ in which the water should not be abstracted.

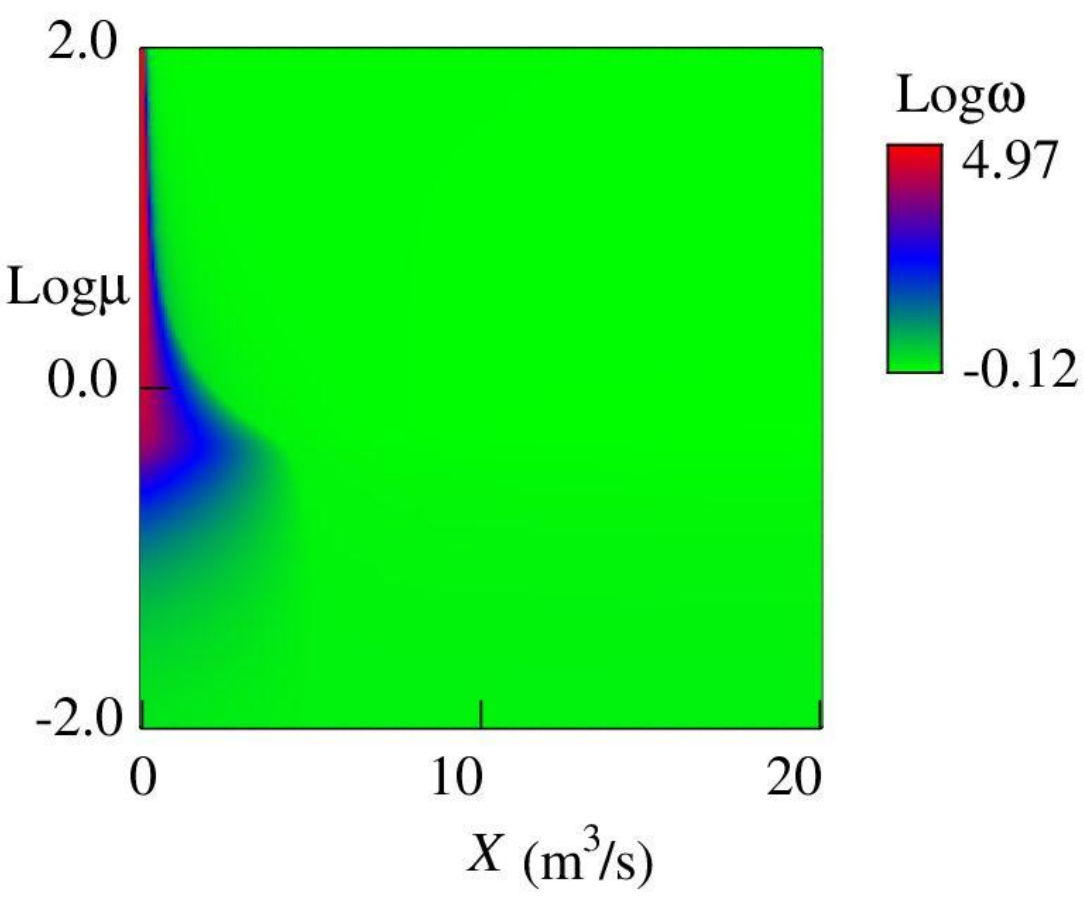

**Figure 9.** The ambiguity under the worst-case optimization as a function of $(X, \mu)$.

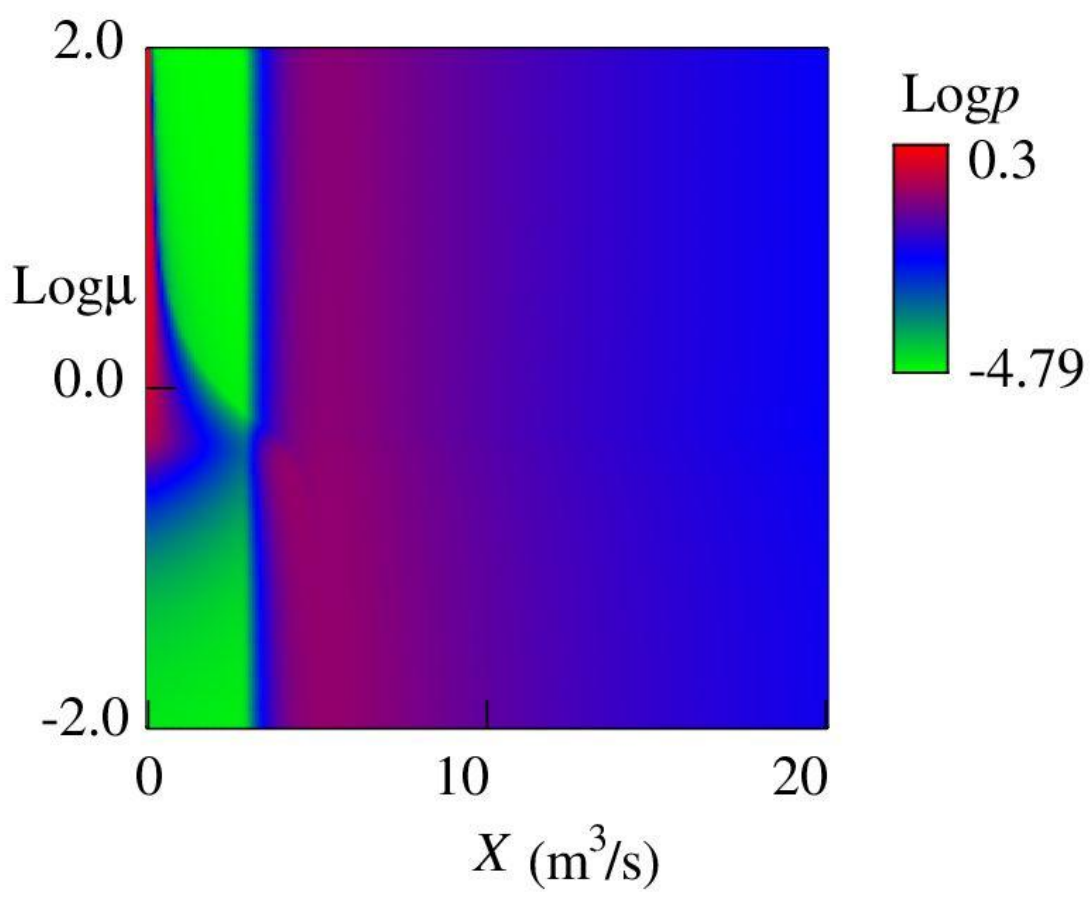

**Figure 10.** The PDF under the worst-case optimization as a function of $(X, \mu)$.

We analyze the similarity between existing diversion ratios and those computed in this case. Razurel et al. (2016) and later studies (Razurel et al., 2018; Perona et al., 2021) utilized the Fermi–Dirac curve to model the diversion ratio. Their parameterized model captures the convex, concave, and sigmoid curves for the diversion ratio as a function of the discharge assuming that it is not large. **Figure 11** plots several optimized diversion ratios for a small discharge with $\alpha = 0.20$ and $\mu = 1.0$. Our model generates the optimal $c = c(X)$ as an almost concave (especially continuous for $X > 0$) and increasing function. This can be classified as a "Fermi standard" curve according to Razurel et al. (2016). We did not find "Fermi inverse" curves or convex and sigmoid "Fermi standard" curves in the literature under the computational conditions for this case. Our computational results imply that the parameterized model of Razurel et al. (2016) may potentially serve as a simple parameterized model for the diversion ratio.

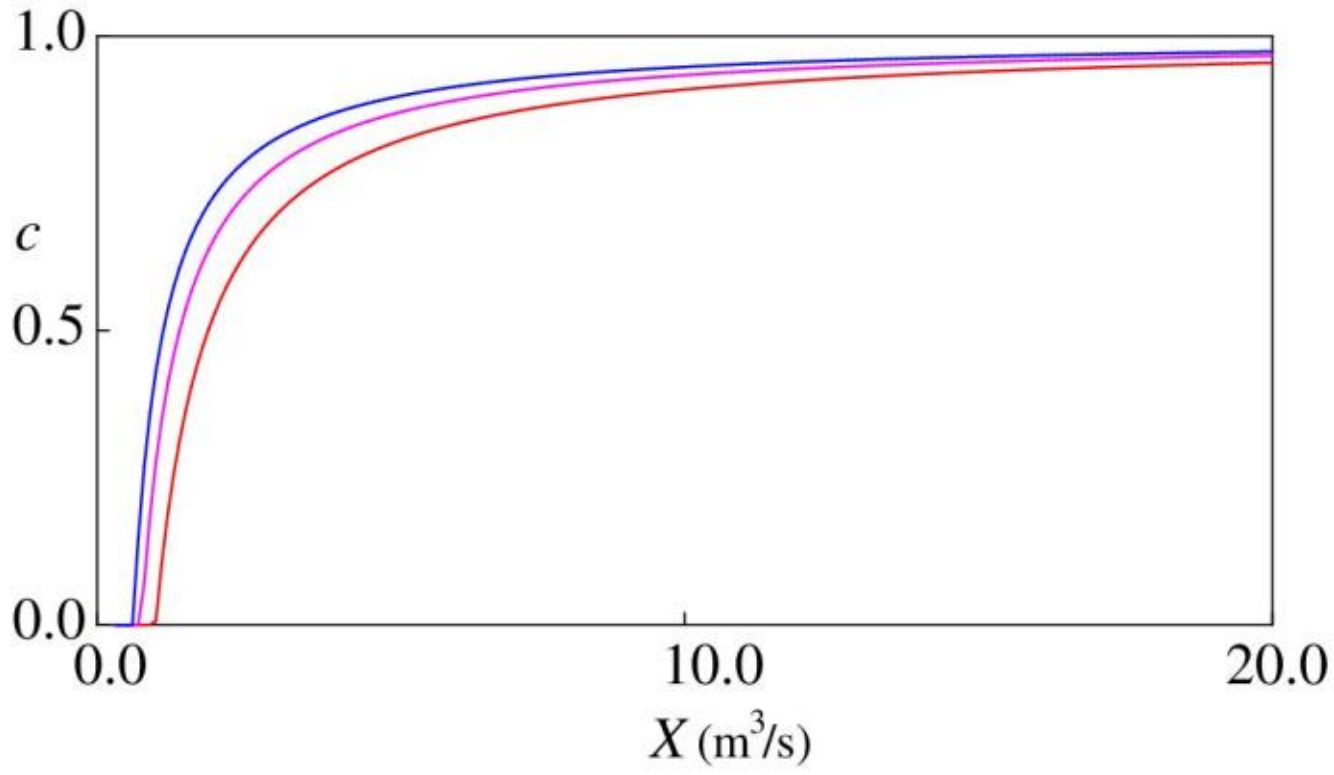

**Figure 11.** The diversion ratio under the worst-case optimization as a function of $X > 0.2$ (m$^3$/s) (Red: $\lambda = 2.5$, Magenta: $\lambda = 3.5$, Blue: $\lambda = 4.5$).

Finally, we analyze the convergence of **Algorithm 1** for different values of $\gamma$. **Figure 12** shows the convergence histories for Case 1 with parameter values for $\gamma = 0.0, 0.1, 0.2, ..., 0.9$. The convergence speed increased with the inertia (i.e., $\gamma$). Hence, adding the inertia to gradient descent was effective at numerically solving the optimization problem. Another key finding is that the convergence history became more unstable as the inertia increases. Therefore, $\gamma$ should not be excessively large. Furthermore, **Algorithm 1** fail when $\gamma \geq 1$, which showed that increasing the inertia too much destroyed the algorithm. Comparable results are obtained for larger $\mu$ (**Figure 13**).

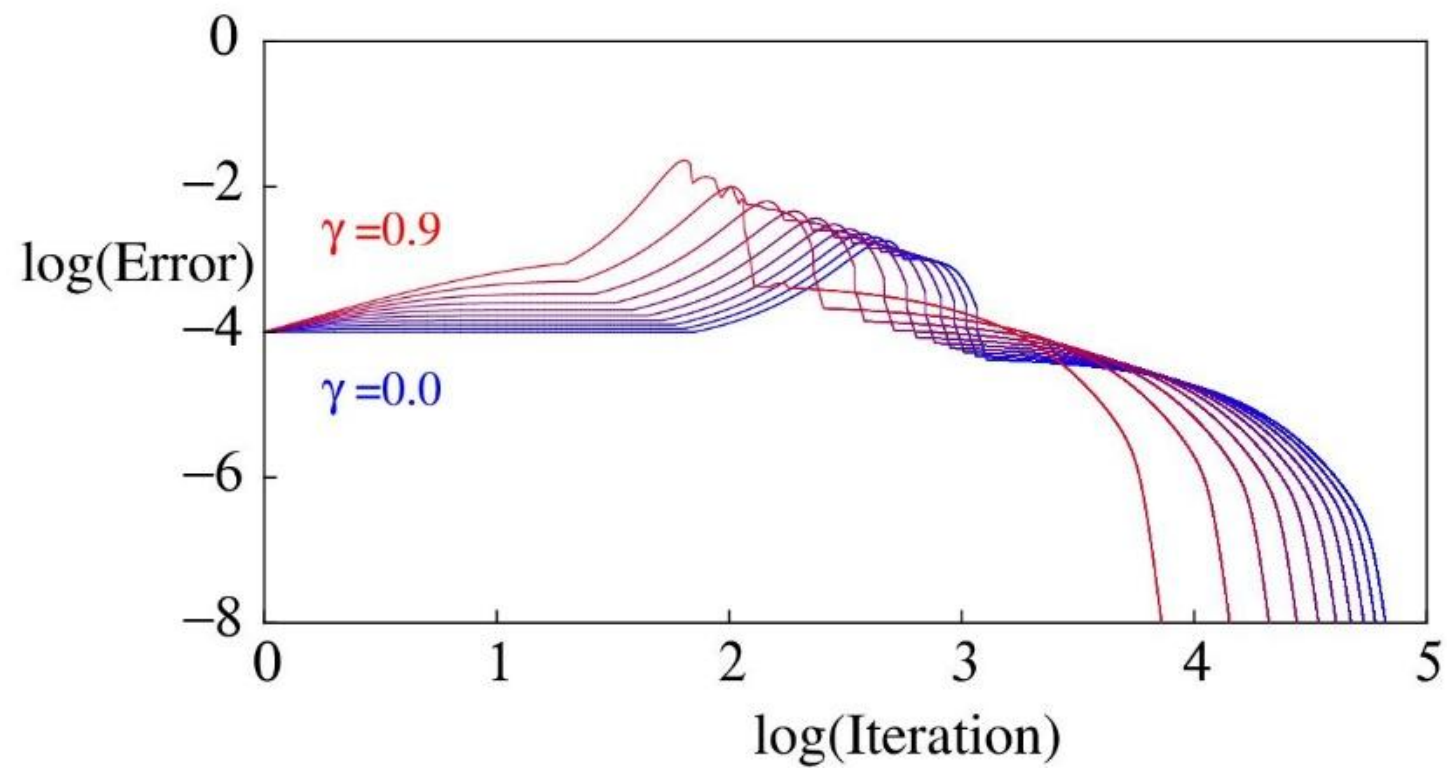

**Figure 12.** Convergence histories of Case 1 for $\gamma = 0.0, 0.1, 0.2, ..., 0.9$ ($\mu = 1$).

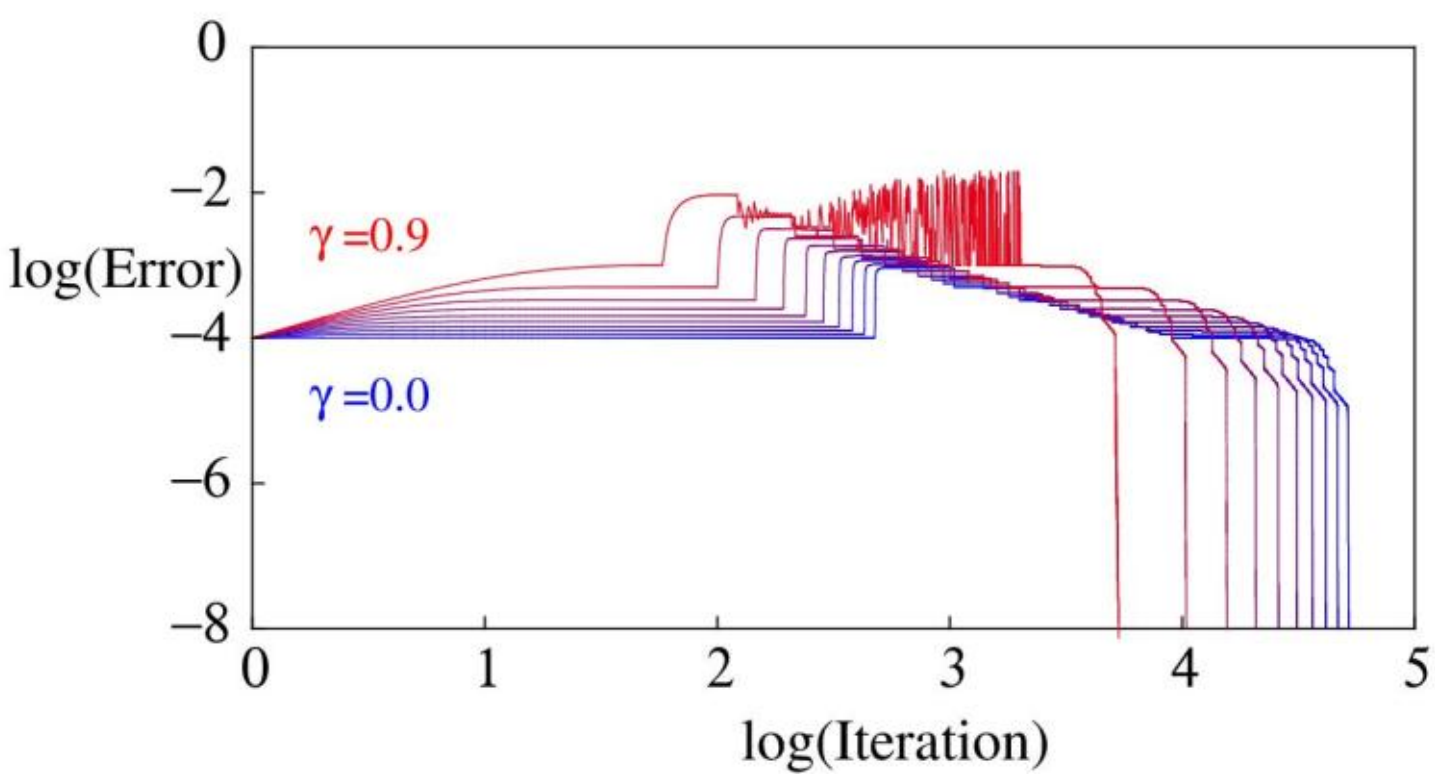

**Figure 13.** Convergence histories of Case 1 for $\gamma = 0.0, 0.1, 0.2, ..., 0.9$ ($\mu = 100$).

### 4.3 Case 2: Control of the high flow

#### 4.3.1 Study site

The operation of flood diversion facilities is a pivotal issue for disaster mitigation of urban areas around a river channel vulnerable to inundation. The Hii River diversion system has been operating since 2013, and it is a huge infrastructure that diverts severe floods from the main channel of the Hii River to the nearby Kando River to mitigate flood damage to downstream urban areas where the ground level is lower than the riverbed (Okada et al., 2013) (**Figure 14**). The brackish Lake Shinji may also trigger prolonged inundation of the surrounding urban area, which is on a plain. The diversion system is activated once or twice a year when the main channel of the Hii River suffers from severe flooding. The discharge data are available at Kamishima station of MLIT, which is 7.7 km upstream of the diversion channel. A lodging weir is equipped at the upstream-end of the diversion channel, which is designed to lodge if the discharge of the main channel exceeds 500 $m^3$/s and rolls at a discharge of 400 $m^3$/s during flood events. Gotoh et al. (2017) reported that the system diverted almost half of the peak flow of a flood event in September 2015. An exactly solvable model was previously proposed for operating the diversion channel under model ambiguity (Yoshioka and Yoshioka, 2022b), but this model assumed that the diversion ratio $c$ is a positive constant. Our proposed model does not *a priori* limit the functional shape, so it is expected to generate more flexible diversion policies that balance the inundation risk of urban areas and diversion efforts under model ambiguity. A major concern for the diversion channel is excessive sediment transport to the Kando River (Gotoh et al., 2017), which would increase with the diversion ratio. **Appendix B** shows that the streamflow discharge at this study site has a long memory, again supporting the use of the MMA process.

**Figure 15** compares the observed and modeled PDFs of the discharge in ordinary and logarithmic scales. Both PDFs are unimodal, and their tails agree well. However, the modeled PDF slightly overestimates the observed PDF near the mode and at large discharge values. The agreement between the tails of the PDFs is important in this case because the focus is on flood diversion. The ambiguity in this context would be introduced by climate change, which may affect the occurrence of extreme rainfall events leading to severe flooding of the main channel.

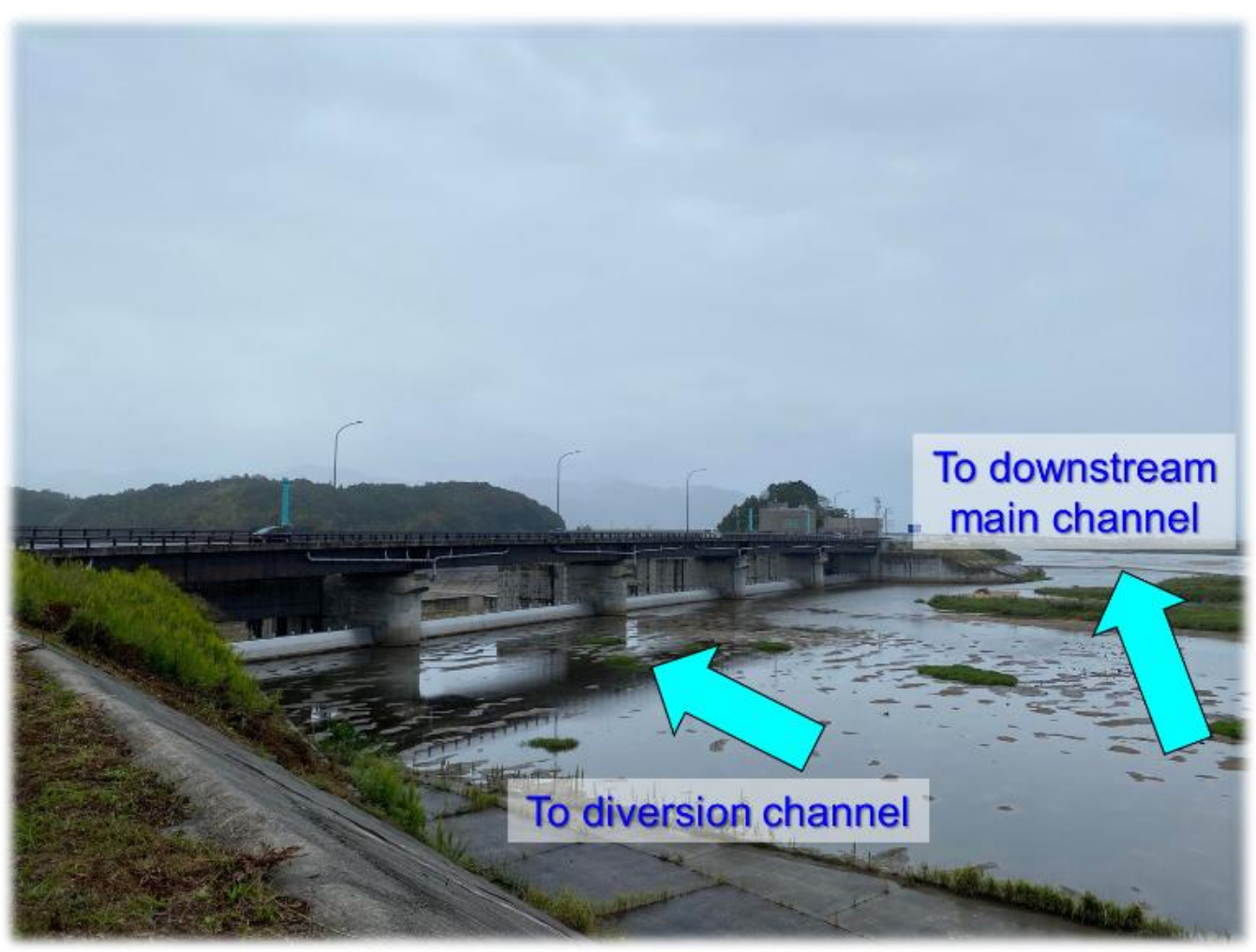

**Figure 14.** The diversion channel taken from upstream. Photo taken by Hashiguchi A. on October 5, 2022.

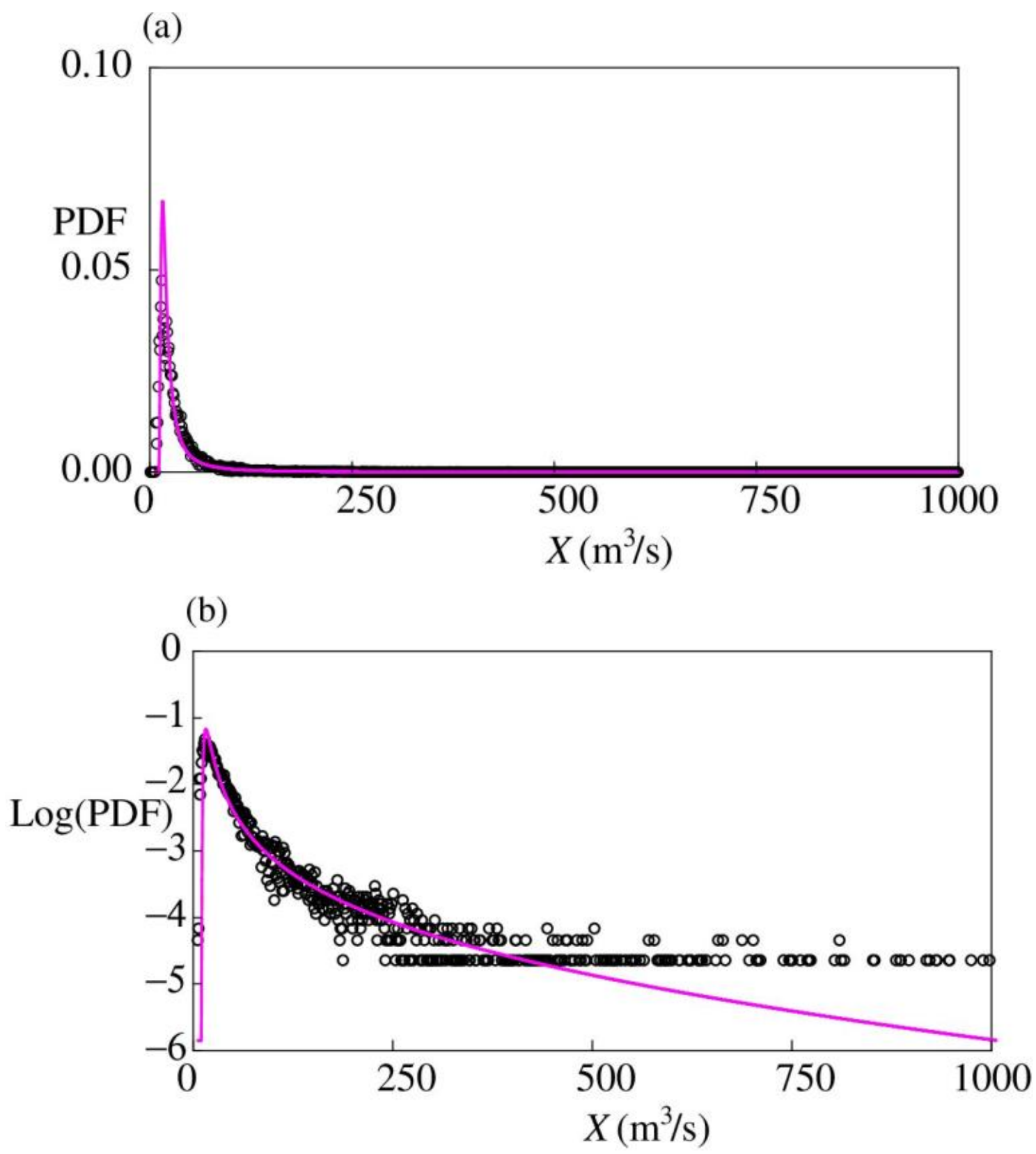

**Figure 15**. Observed and modelled PDFs on **(a)** ordinary scale and **(b)** logarithmic scale.

#### 4.3.2 Computation and sensitivity analysis

We evaluated the impact of the parameter $\lambda$, which controls the CVaR term to penalize the flood risk downstream of the main channel. **Figures 16** and **17** show the computed diversion ratio $c$ under the optimality for the quantiles $\beta$ of 0.99 and 0.999, respectively. We set $\mu = 10$. The computational results suggest that treating the flood risk more seriously by choosing a larger $\lambda$ would result in a more conservative diversion policy with larger $c$. Interestingly, we found that the optimal $c = c(X)$ as a function of the discharge $X$ seems to be continuous, but it is not continuously differentiable at few points. For a relatively small $\lambda$ with $\beta = 0.99$, it is piecewise-linear when it is positive and becomes smoother as $\lambda$ increases. This qualitative transition is not observed for $\beta = 0.999$ under the specified computational conditions. In all cases, $c = c(X)$ seems nondifferentiable at the point where it transitioned from zero to positive. This indicates that streamflow management should only activate flood diversion when the discharge is greater than some threshold value determined from the hydrological properties of the river flow and preferences of the decision-maker.

We analyzed the impact of the parameter $\mu$ on the tail of the discharge PDF by setting $\lambda = 0.0001$. **Figure 18** shows the computed ambiguity under the worst-case optimization for different values of $\mu$ under $\beta = 0.999$. **Figure 19** shows the corresponding computed PDF. Increasing the ambiguity aversion by taking a smaller $\mu$ results in a fatter tail of the PDF. Parts of the PDF far from the tail are not affected by the model ambiguity because they did not critically affect the flood risk. Thus, the unimodal shape of the PDF is not broken in the worst case. Finally, we analyzed the statistical performance of the optimized diversion ratio by using a Monte Carlo simulation based on the Markovian lift. The algorithm is presented in **Appendix C**. We focus on the maximum diversion ratio and maximum diverted discharge within 1 year, both of which can be used as design criteria for flood diversion facilities. We set $\beta = 0.999$ and three cases for the weight $\lambda$ as examples: 0.000001, 0.0000008, and 0.000016. The corresponding diversion ratios are shown in **Figure 19**. The Monte Carlo simulation was run for 1,000 years, and the maximum diverted discharge $cX$ was obtained for each year.

**Figure 20** shows a sample path of the simulated discharge in the first year and the corresponding paths of diverted discharges for different values of $\lambda$. The floods are diverted accordingly to different magnitudes in the different cases. **Figure 21** shows the histograms of the yearly maximum diverted discharges for different values of $\lambda$. **Table 1** presents their statistics. The yearly maximum diverted discharges follow highly non-Gaussian distributions with mild local maxima and large boundary peaks. Their probability densities were trimodal for relatively small $\lambda$, while the peak of the lower part of the probability density disappears with increasing $\lambda$. However, the peak of the lower part was affected by small diversion ratios (<0.00001), which may be attributed to regularization. Thus, the computational results suggest that designing and operating a flood diversion system is critically affected by the ambiguity aversion of the decision-maker. In particular, the ability of the system should be designed to accommodate expected realizations of the maximum diverted discharge.

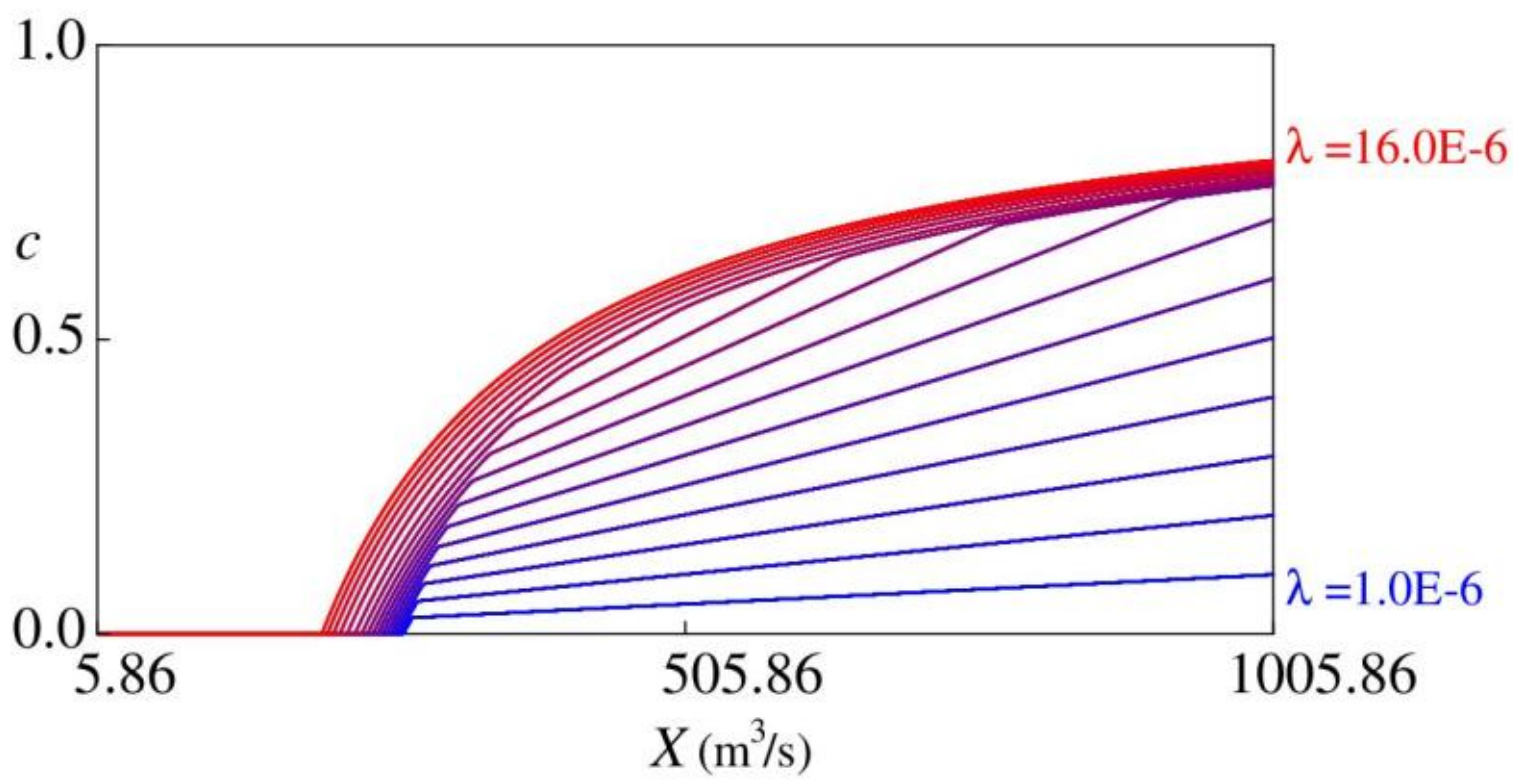

**Figure 16**. Computed diversion ratio under the worst-case optimization with $\beta = 0.99$. The plots are based on $\lambda = 0.000001$ to $\lambda = 0.000016$ with the increment 0.000001 from the bottom to the top.

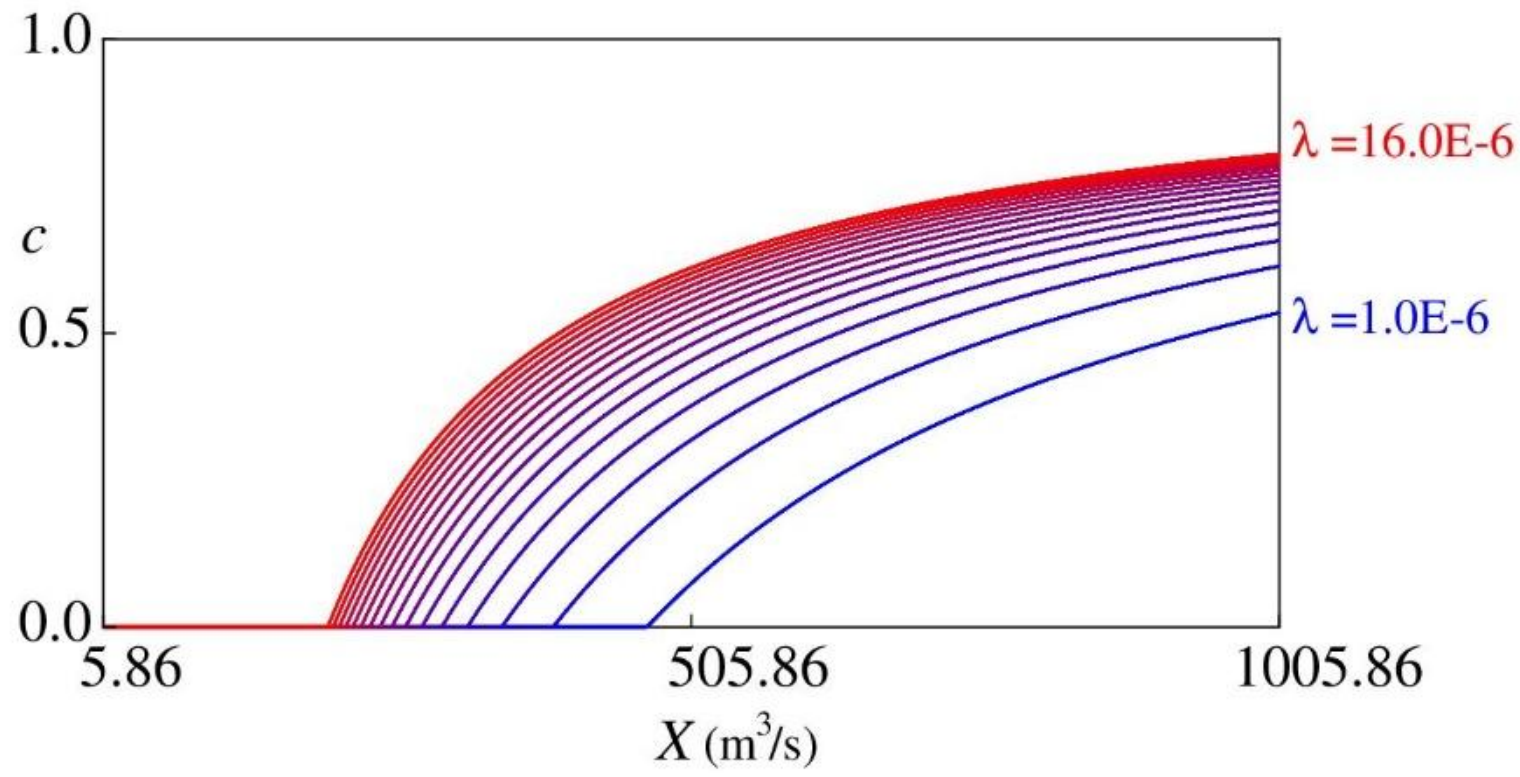

**Figure 17**. Computed diversion ratio under the worst-case optimization with $\beta = 0.999$ . The plots are based on $\lambda = 0.000001$ to $\lambda = 0.000016$ with the increment $0.000001$ from the bottom to the top.

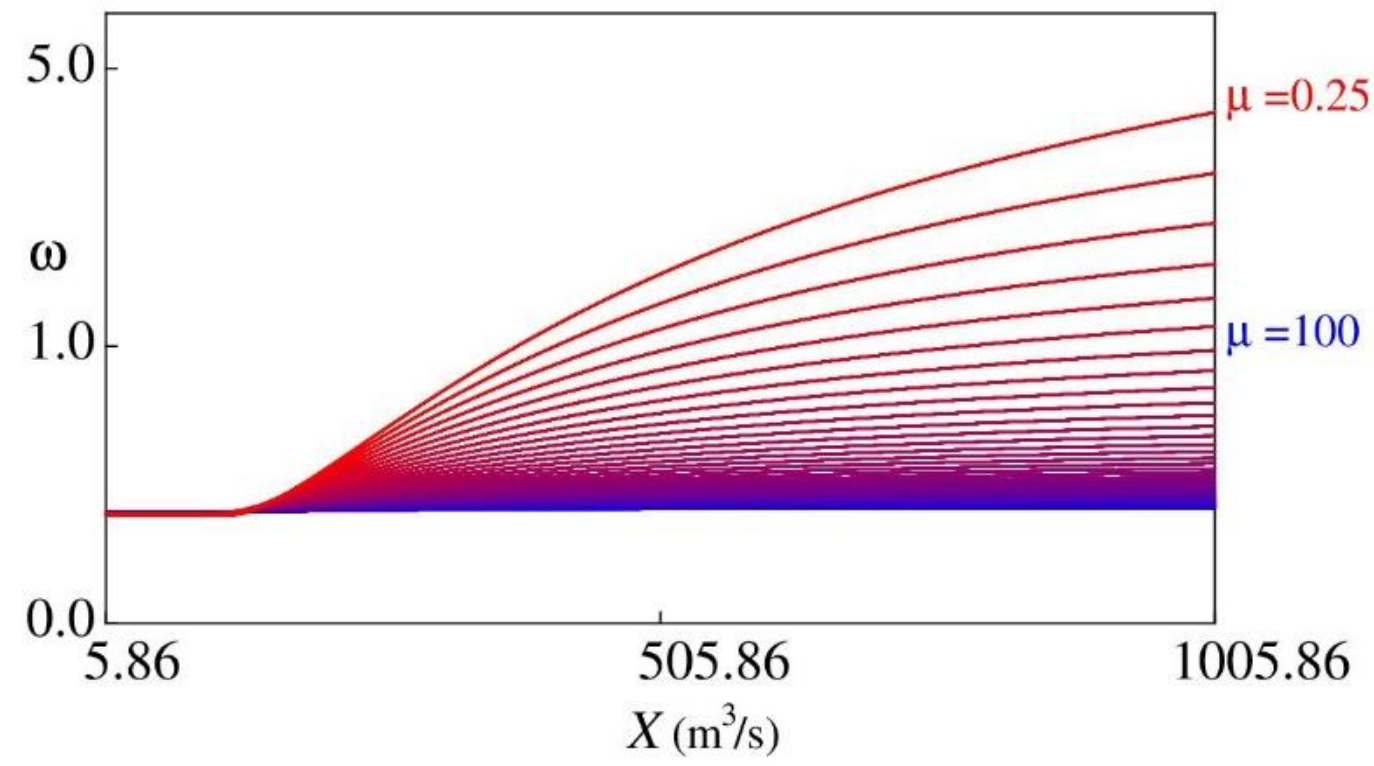

**Figure 18**. Computed ambiguity under the worst-case optimization with $\beta = 0.999$ . The plots are based on $\mu = 10^{1-i/25}$ ( $i = 0,1,2,...,40$ ) from the bottom to the top.

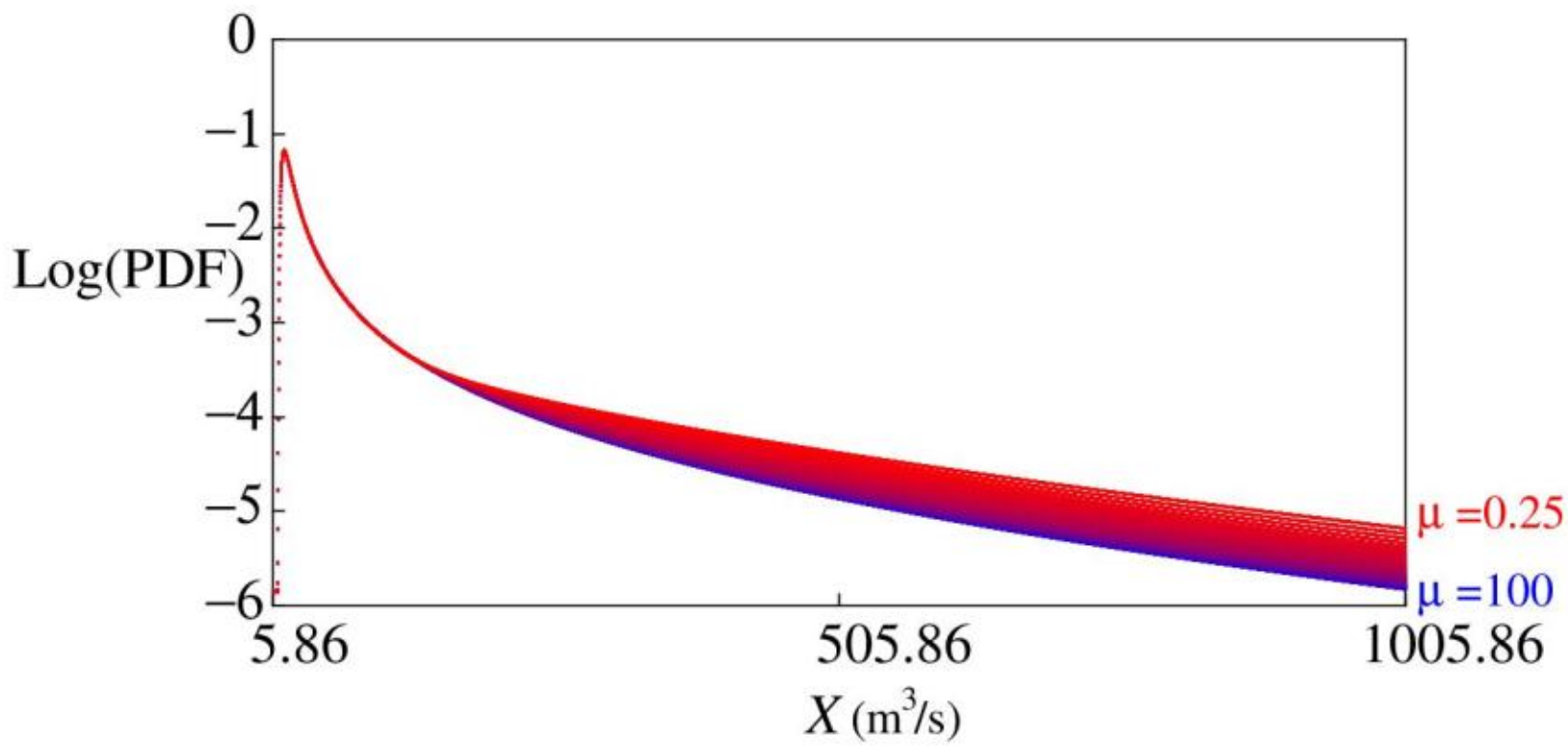

**Figure 19**. Computed PDF under the worst-case optimization with $\beta = 0.999$ . The plots are based on $\mu = 10^{1-i/25}$ ( $i = 0,1,2,...,40$ ) from the bottom to the top.

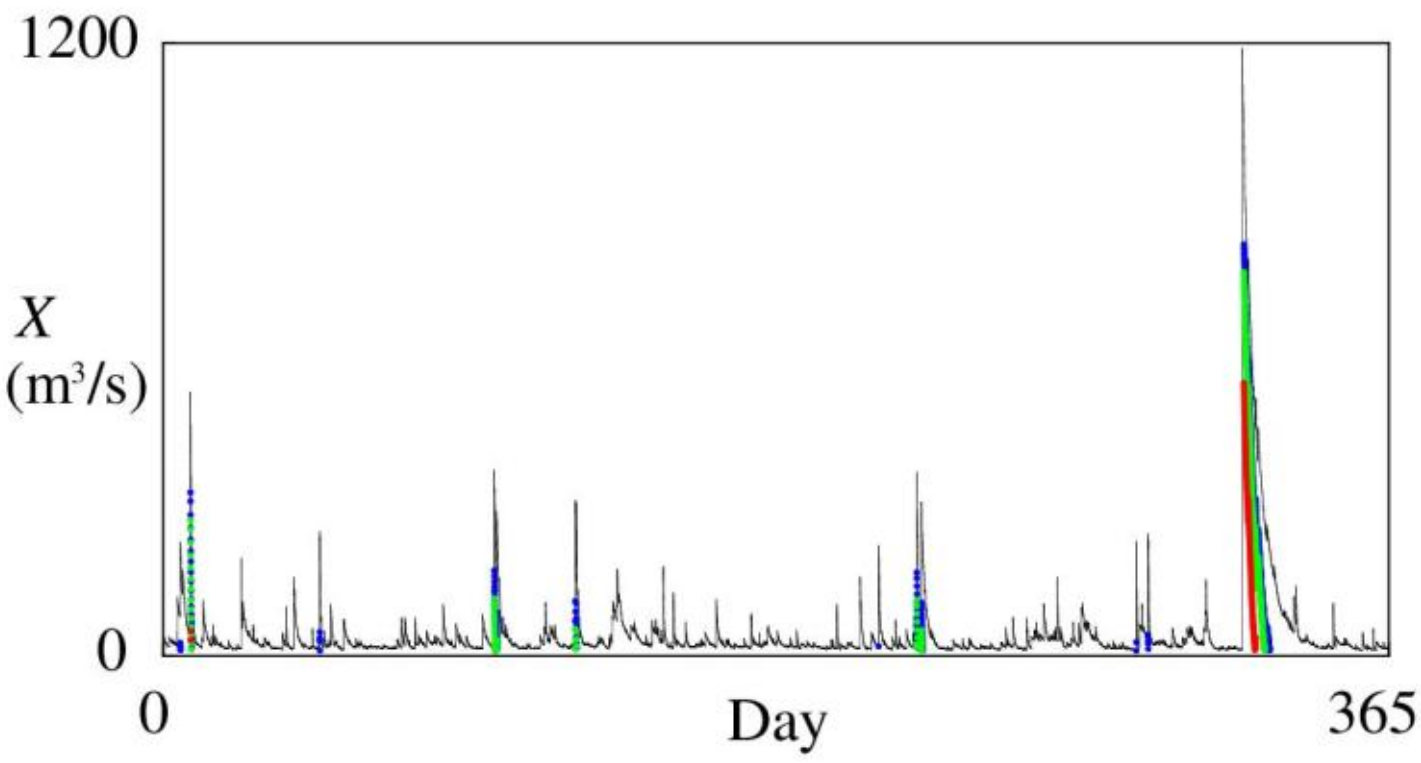

**Figure 20.** The streamflow discharge and the corresponding diverged discharges (Red: $\lambda = 0.000001$, Green: $\lambda = 0.000008$, Blue: $\lambda = 0.000016$).

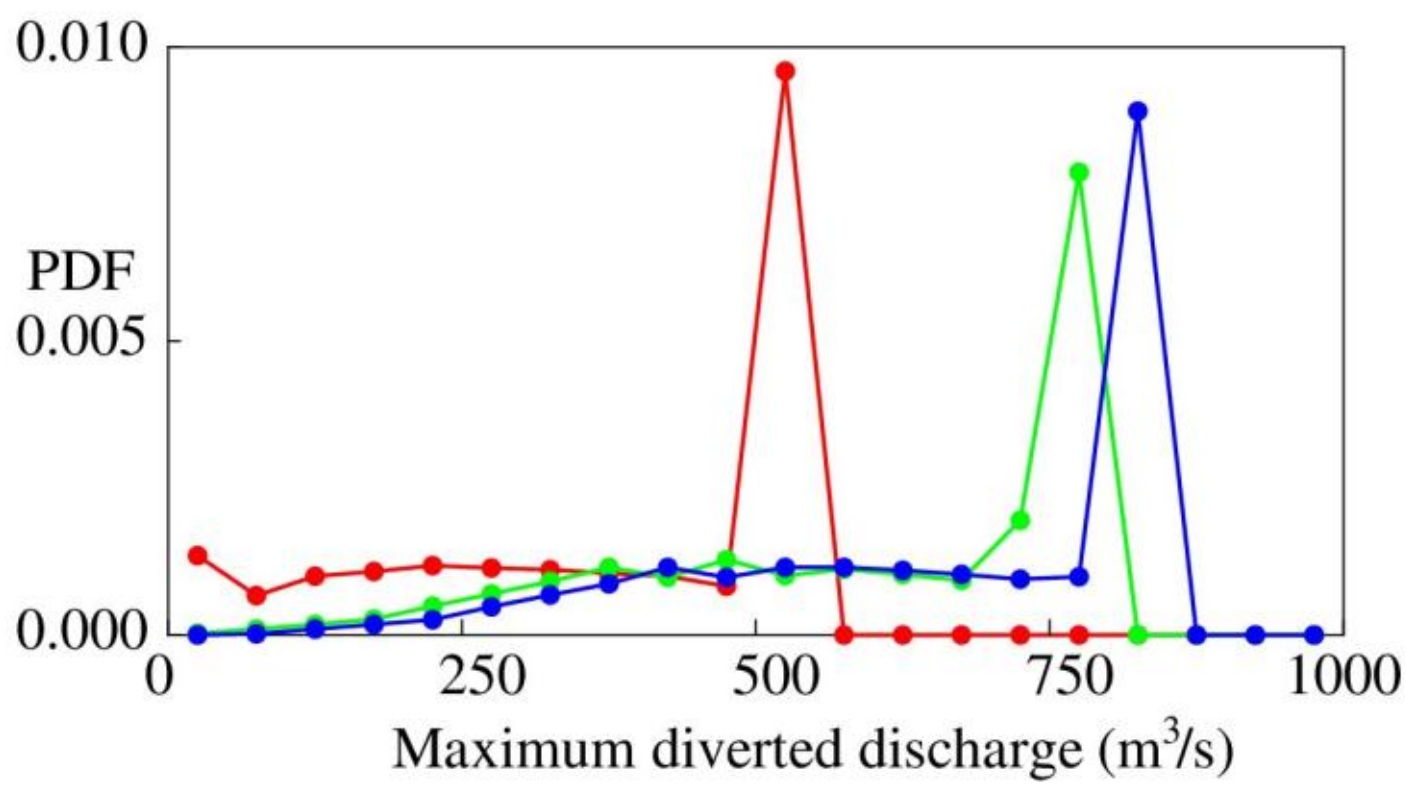

**Figure 21**. Computed PDFs of the yearly maximum diverted discharges (Red: $\lambda = 0.000001$, Green: $\lambda = 0.000008$, Blue: $\lambda = 0.000016$).

**Table 1**. Statistics of the computed yearly maximum diverted discharges.

| | $\lambda = 0.000001$ | $\lambda = 0.000008$ | $\lambda = 0.000016$ |
|---|---|---|---|
| Average (m$^3$/s) | 3.83.E+02 | 5.98.E+02 | 6.52.E+02 |
| Standard deviation (m$^3$/s) | 1.79.E+02 | 1.84.E+02 | 1.84.E+02 |
| Skewness (-) | 3.19.E+04 | 3.40.E+04 | 3.40.E+04 |
| Kurtosis (-) | −7.71.E-01 | −8.81.E-01 | −8.81.E-01 |

## 5. Conclusion

We presented a consistent modeling framework for distributionally robust optimization by covering the underlying MMA process, the Markovian lift, the conversion of the characteristic function to a PDF, the formulation of the optimization problem, and the algorithm to find the numerical solution. A key to the modeling framework is the affine MMA process, which provides a flexible mathematical model with a closed-form characteristic function from which the PDF can be obtained. Another key is the strict convexity of the regularized optimized problem, which obtains unique values for the optimal decision variables as well as the worst-case model ambiguity. The supCBI processes were successfully identified at two sites with mutually different objectives, and the optimization problem was computed at each site to obtain the optimal water diversion policies. In this paper, we report that the proposed method is applicable to problems with models identified from data and explicitly quantify how the optimal diversion ratios and controlled streamflow discharge depend on the parameter values in the models.

The convexity of the optimization problem is preserved if the CVaR is replaced by a spectral risk measure (Rockafellar and Royset, 2018; Guo and Xu, 2022) or adjusted entropic value-at-risk (Zou et al., 2022). Model ambiguity can be evaluated by a generalized method, such as the Tsallis divergence (Fu et al., 2022) with a suitable modification, further highlighting the extendibility of our framework. The modeling framework and optimization of affine processes can be extended to Markov-modulated cases where the coefficients of SDEs depend exogenously on another SDE of an affine or polynomial type (Kurt and Frey, 2022) if the characteristic function is available in a closed form. Multistage problems (Dowson et al., 2022) will be studied in the future based on the recursive use of the single-stage problem formulated in this paper. Our mathematical and computational methods would be applicable to the optimization of other indices arising in different research areas, including finance and economics (e.g., time series of asset prices) and civil engineering (e.g., time series of pollutant loads), where MMA and related stochastic processes play a role. Properties of the proposed optimization model can be analyzed more in detail if the notion of Orlicz space (Part IV in Rubshtein et al. (2016)) is applied to each expectation in the objective function. Such an application with a different viewpoint from the proposed one is found in Yoshioka and Yoshioka (2024).

## Appendices

### Appendix A: Proofs of propositions and a lemma

***Proof of Proposition 2***

Firstly, the right-hand side of (10) exists. Indeed, by (9) we have

$$\begin{aligned}
&\frac{\partial}{\partial \xi}\exp\left(A\int_0^{+\infty}\frac{\pi(\mathrm{d}\rho)}{\rho}\int_0^{+\infty}\int_0^{+\infty}\left(\exp\left(\xi z e^{-t}\right)-1\right)v(\mathrm{d}z)\mathrm{d}t\right)\\
&= A'M_{X_t}(\xi)\frac{\partial}{\partial \xi}\int_0^{+\infty}\left(1-\left(1-\frac{\xi e^{-t}}{\beta_v}\right)^{\alpha_v}\right)\mathrm{d}t\\
&= -A'M_{X_t}(\xi)\int_0^{+\infty}\frac{\partial}{\partial \xi}\left(1-\frac{\xi e^{-t}}{\beta_v}\right)^{\alpha_v}\mathrm{d}t\\
&= -\alpha_v A'M_{X_t}(\xi)\int_0^{+\infty}\left(1-\frac{\xi e^{-t}}{\beta_v}\right)^{\alpha_v-1}\frac{e^{-t}}{\beta_v}\mathrm{d}t
\end{aligned} \quad . \tag{58}$$

Here, we used the boundedness of the integrand between the second and third lines. The last integral of (58) is found analytically, thereby

$$\begin{aligned}
-\alpha_v A'M_{X_t}(\xi)\int_0^{+\infty}\left(1-\frac{\xi e^{-t}}{\beta_v}\right)^{\alpha_v-1}\frac{e^{-t}}{\beta_v}\mathrm{d}t &= \alpha_v A'M_{X_t}(\xi)\left[\frac{1}{\alpha_v\xi}\left(1-\frac{\xi e^{-t}}{\beta_v}\right)^{\alpha_v}\right]_0^{+\infty}\\
&= A'M_{X_t}(\xi)\left[\frac{1}{\xi}\left(1-\frac{\xi e^{-t}}{\beta_v}\right)^{\alpha_v}\right]_0^{+\infty}\\
&= A'M_{X_t}(\xi)\frac{1}{\xi}\left(1-\left(1-\frac{\xi}{\beta_v}\right)^{\alpha_v}\right)
\end{aligned} \quad . \tag{59}$$

We have

$$\lim_{\xi\to+0}\frac{1}{\xi}\left(1-\left(1-\frac{\xi}{\beta_v}\right)^{\alpha_v}\right)=\frac{\alpha_v}{\beta_v}>0 \tag{60}$$

by the Taylor expansion

$$1-\left(1-\frac{\xi}{\beta_v}\right)^{\alpha_v}=\frac{\alpha_v}{\beta_v}\xi+O\left(\xi^2\right) \quad \text{for small } \xi>0 . \tag{61}$$

Therefore, $M_{X_t}(\xi)$ and hence the last line of (59) is bounded if $\xi \le \beta_v$. Then, we can exchange the partial differentiation and expectation:

$$\frac{\partial}{\partial \xi}\mathbb{E}\left[e^{\xi X_t}\right]=\mathbb{E}\left[\frac{\partial}{\partial \xi}e^{\xi X_t}\right]=\mathbb{E}\left[X_t e^{\xi X_t}\right], \tag{62}$$

yielding (10). For (11), it is sufficient to check $\lim_{\xi\to+0}\frac{\mathrm{d}}{\mathrm{d}\xi}\frac{1}{\xi}\left(1-\left(1-\frac{\xi}{\beta_v}\right)^{\alpha_v}\right)<+\infty$. This follows from (61).

□

***Proof of Proposition 3***

Firstly, we show the sign inequalities

$$F(\omega)\equiv-\omega+B\int_0^{+\infty}\left(\exp(\omega z)-1\right)v(\mathrm{d}z)\le 0 \quad \text{and} \quad \frac{\mathrm{d}F(\omega)}{\mathrm{d}\omega}\le 0 \tag{63}$$

for $\omega \le \beta_c$. Indeed,

$$F(\omega) = -\omega + \frac{\beta_v}{\alpha_v} BM_1 \left(1 - \left(1 - \frac{\omega}{\beta_v}\right)^{\alpha_v}\right) \text{ and } \frac{\mathrm{d}F(\omega)}{\mathrm{d}\omega} = -1 + BM_1 \left(1 - \frac{\omega}{\beta_v}\right)^{\alpha_v - 1}, \tag{64}$$

showing that $\frac{\mathrm{d}F(\omega)}{\mathrm{d}\omega} \le 0$ for $0 \le \omega \le \beta_c$, especially $\frac{\mathrm{d}F(\omega)}{\mathrm{d}\omega} < 0$ for $0 \le \omega < \beta_c$. Further, $F(0) = 0$ and hence $F(\omega) \le 0$ for $\omega \le \beta_c$. Due to the smoothness (Lipschitz continuity) and boundedness of $F(\omega)$ for $\omega \le \beta_c$, the classical unique existence and comparison theorems show that the initial value problem (17) admits a unique global solution $\varpi$ strictly smaller than $\beta_v$. Using this $\varpi$, the moment generating function, if it exists, is expressed as (16).

The existence of the right-hand side of (16) holds true if $\int_0^{+\infty}\int_0^{+\infty}\left(\exp\left(\varpi(t)z\right)-1\right)v(\mathrm{d}z)\mathrm{d}t < +\infty$. This follows if we can find the bound of the form

$$\varpi(t) \le \beta_c \exp(-\sigma t), \quad t \ge 0 \tag{65}$$

for a suitable constant $\sigma > 0$. Such a $\sigma$ exists as proven below. Set

$$G(\omega) \equiv -\sigma\omega - F(\omega) = (1-\sigma)\omega - \frac{\beta_v}{\alpha_v} BM_1 \left(1 - \left(1 - \frac{\omega}{\beta_v}\right)^{\alpha_v}\right), \quad \omega \le \beta_c. \tag{66}$$

We then obtain

$$\frac{\mathrm{d}G(\omega)}{\mathrm{d}\omega} = (1-\sigma) - BM_1 \left(1 - \frac{\omega}{\beta_v}\right)^{\alpha_v - 1}, \quad \omega \le \beta_c, \tag{67}$$

showing that $\frac{\mathrm{d}G(\omega)}{\mathrm{d}\omega} \le 0$ for $\omega \le \beta_c$ if $\frac{\mathrm{d}G(\beta_c)}{\mathrm{d}\omega} \le 0$, but the last condition is satisfied since

$$\begin{aligned}
\frac{\mathrm{d}G(\beta_c)}{\mathrm{d}\omega} &= 1 - \sigma - BM_1 \left(1 - \frac{\beta_v\left(1 - (BM_1)^{\frac{1}{1-\alpha_v}}\right)}{\beta_v}\right)^{\alpha_v - 1} \\
&= 1 - \sigma - BM_1 \left(1 - \left(1 - (BM_1)^{\frac{1}{1-\alpha_v}}\right)\right)^{\alpha_v - 1} \\
&= 1 - \sigma - BM_1 \left((BM_1)^{\frac{1}{1-\alpha_v}}\right)^{\alpha_v - 1} \\
&= -\sigma < 0
\end{aligned}. \tag{68}$$

Therefore, we obtain $G(\omega) \ge 0$ for $\omega \le \beta_c$ if $G(\beta_c) \ge 0$ because $G$ is non-increasing for $\omega \le \beta_c$; namely, if

$$
\begin{aligned}
G(\beta_c) &= (1-\sigma)\beta_c - \frac{\beta_v}{\alpha_v} BM_1 \left( 1 - \left( 1 - \frac{\beta_v \left( 1 - (BM_1)^{\frac{1}{1-\alpha_v}} \right)}{\beta_v} \right)^{\alpha_v} \right) \\
&= (1-\sigma)\beta_c - \frac{\beta_v}{\alpha_v} BM_1 \left( 1 - \left( 1 - \left( 1 - (BM_1)^{\frac{1}{1-\alpha_v}} \right) \right)^{\alpha_v} \right) \\
&= (1-\sigma)\beta_c - \frac{\beta_v}{\alpha_v} BM_1 \left( 1 - (BM_1)^{\frac{\alpha_v}{1-\alpha_v}} \right) \\
&= (1-\sigma)\beta_v \left( 1 - (BM_1)^{\frac{1}{1-\alpha_v}} \right) - \frac{\beta_v}{\alpha_v} \left( BM_1 - (BM_1)^{\frac{1}{1-\alpha_v}} \right) \\
&= \beta_v \left[ (1-\sigma)\left( 1 - (BM_1)^{\frac{1}{1-\alpha_v}} \right) - \frac{1}{\alpha_v} \left( BM_1 - (BM_1)^{\frac{1}{1-\alpha_v}} \right) \right] \\
&\geq 0
\end{aligned}
\qquad (69)
$$

By (69), $G(\beta_c) \geq 0$ if we choose $\sigma$ as

$$
0 < \sigma \leq 1 - \frac{1}{\alpha_v} \frac{BM_1 - (BM_1)^{\frac{1}{1-\alpha_v}}}{1 - (BM_1)^{\frac{1}{1-\alpha_v}}}, \qquad (70)
$$

which is possible if $B$ is small. Consequently, there exists a small $\sigma > 0$ such that $G(\omega) \geq 0$ for $\omega \leq \beta_c$ if $B$ is sufficiently small. We obtain $-\sigma\omega \geq F(\omega)$ for $\omega \leq \beta_c$ in this case, and hence the comparison of solutions to the two differential equations with the right-hand sides $F(\omega)$ and $-\sigma\omega$; the solution to the latter is larger than that of the former $\varpi$, yielding (65). Under the same condition, we obtain the desired estimate (65). Then, the moment generating function exist, and further that the latter equals $\mathbb{E}\left[e^{\xi X_t}\right]$.

□

***Proof of Proposition 4***

We follow the proof of Theorem 6 of Liu et al. (2022). Given $u, w \geq 0, c \in \mathfrak{C}$, by the classical Jensen's inequality and the non-negativity of $g$, we have

$$
\begin{aligned}
&-\lambda u + \eta w + H(u,w,c) \\
&= -\lambda u + \eta w + \mu \ln\left( \int e^{\frac{F(x,u,w,c)}{\mu}} p(x)\mathrm{d}x \right) \\
&\geq -\lambda u + \eta w + \int F(x,u,w,c)\, p(x)\mathrm{d}x \\
&\geq -\lambda u + \eta w + \int \left( \frac{\lambda}{\alpha} \max\left\{ u - (1-c(x))x, 0 \right\} + \frac{\eta}{1-\beta} \max\left\{ (1-c(x))x - w, 0 \right\} \right) p(x)\mathrm{d}x \\
&\equiv \underline{H}(u,w,c)
\end{aligned}
\qquad (71)
$$

We then obtain

$$
\begin{aligned}
&\liminf_{w\to+\infty} \underline{H}(u,w,c) \\
&= \liminf_{w\to+\infty}\left(\eta w + \int \frac{\eta}{1-\beta}\max\left\{(1-c(x))x - w, 0\right\} p(x)\mathrm{d}x\right) \\
&+\left(-\lambda u + \int \frac{\lambda}{\alpha}\max\left\{u-(1-c(x))x, 0\right\} p(x)\mathrm{d}x\right) \\
&\geq \eta \liminf_{w\to+\infty} w + \left(-\lambda u + \int \frac{\lambda}{\alpha}\max\left\{u-(1-c(x))x, 0\right\} p(x)\mathrm{d}x\right) \\
&= +\infty
\end{aligned}
\tag{72}
$$

and

$$
\begin{aligned}
&\liminf_{u\to+\infty} \underline{H}(u,w,c) \\
&= \liminf_{u\to+\infty}\left(-\lambda u + \eta w + \int\left(\frac{\lambda}{\alpha}\max\left\{u-(1-c(x))x, 0\right\} + \frac{\eta}{1-\beta}\max\left\{(1-c(x))x - w, 0\right\}\right) p(x)\mathrm{d}x\right) \\
&\geq \liminf_{u\to+\infty}\left(-\lambda u + \int \frac{\lambda}{\alpha}\max\left\{u-(1-c(x))x, 0\right\} p(x)\mathrm{d}x\right) \\
&= \frac{\lambda}{\alpha}\liminf_{u\to+\infty}\left(-\alpha u + \int \max\left\{u-(1-c(x))x, 0\right\} p(x)\mathrm{d}x\right) \\
&\geq \frac{\lambda}{\alpha}\liminf_{u\to+\infty}\left(-\alpha u + \int \left(u-(1-c(x))x\right) p(x)\mathrm{d}x\right) \\
&= \frac{\lambda}{\alpha}\liminf_{u\to+\infty}\left((1-\alpha)u - \int (1-c(x)) x p(x)\mathrm{d}x\right) \\
&= +\infty
\end{aligned}
\quad . \tag{73}
$$

Therefore, both separately and simultaneously sending $u, v \to +\infty$ result in $W \to +\infty$. Further, we have

$$
W \leq \underline{H}(0,0,c\equiv 1) = \mu \ln \mathbb{E}\left[e^{\frac{F(X,0,0,1)}{\mu}}\right] = \mu \ln \mathbb{E}\left[e^{\frac{1}{\mu}g(|1-\hat{c}|)}\right] = g(|1-\hat{c}|) < +\infty \,, \tag{74}
$$

showing that the optimal decision variables must be contained in a compact set. Consequently, we can replace the admissible set of decision variables from $u, w \geq 0, c \in \mathfrak{C}$ to $(u,w) \in D, c \in \mathfrak{C}$ with some compact set $D \subset [0,+\infty)^2$ without changing $W$.

□

***Proof of Proposition 5***

The lower semi-continuity of $H = H(u,w,c)$ with respect to $c \in \mathfrak{C}$ for each given $(u,w) \in [0,+\infty)^2$ implies (Section 2.6 of Galewski (2021))

$$
H(u,w,c) \leq \liminf_{n\to+\infty} H(u,w,c_n) \tag{75}
$$

for each sequence $c_n \in \mathfrak{C}$ ( $n = 1,2,3,...$ ) converging to $c$. We prove this inequality. We use the elementary inequality

$$
|\ln x_1 - \ln x_2| \leq |x_1 - x_2|, \quad x_1, x_2 \geq 1 \,. \tag{76}
$$

Given $c \in \mathfrak{C}$ and its converging sequence $c_n \in \mathfrak{C}$ ( $n = 1,2,3,...$ ), by $F \geq 0$ we have

$$\begin{aligned}
H(u,w,c) &= \mu \ln\left( \int e^{\frac{F(x,u,w,c)}{\mu}} p(x)\mathrm{d}x \right) \\
&\le \left| \mu \ln\left( \int e^{\frac{F(x,u,w,c)}{\mu}} p(x)\mathrm{d}x \right) - \mu \ln\left( \int e^{\frac{F(x,u,w,c_n)}{\mu}} p(x)\mathrm{d}x \right) \right| + H(u,w,c_n) \\
&\le \left| \mu \int e^{\frac{F(x,u,w,c)}{\mu}} p(x)\mathrm{d}x - \mu \int e^{\frac{F(x,u,w,c_n)}{\mu}} p(x)\mathrm{d}x \right| + H(u,w,c_n) \\
&\le \mu \int \left| e^{\frac{F(x,u,w,c)}{\mu}} - e^{\frac{F(x,u,w,c_n)}{\mu}} \right| p(x)\mathrm{d}x + H(u,w,c_n)
\end{aligned} \quad . \tag{77}$$

Applying the Lipschitz continuity (38) of $F$ and a Cauchy–Schwarz inequality along with **Assumption 1** yields

$$\begin{aligned}
H(u,w,c) &\le \mu \int \left| e^{\frac{F(x,u,w,c)}{\mu}} - e^{\frac{F(x,u,w,c_n)}{\mu}} \right| p(x)\mathrm{d}x + H(u,w,c_n) \\
&\le \mu \int e^{\frac{1}{\mu}\left( g(1) + \frac{\lambda U}{\alpha} + \frac{\eta}{1-\beta} X \right)} \left| \frac{F(x,u,w,c)}{\mu} - \frac{F(x,u,w,c_n)}{\mu} \right| p(x)\mathrm{d}x + H(u,w,c_n) \\
&= e^{\frac{1}{\mu}\left( g(1) + \frac{\lambda U}{\alpha} \right)} \int e^{\frac{1}{\mu}\frac{\eta}{1-\beta} X} \left| F(x,u,w,c) - F(x,u,w,c_n) \right| p(x)\mathrm{d}x + H(u,w,c_n) \\
&\le L' \int e^{\frac{1}{\mu}\frac{\eta}{1-\beta} X} (x + C') \left| c(x) - c_n(x) \right| p(x)\mathrm{d}x + H(u,w,c_n) \\
&\le L' \sqrt{\int (x + C')^2 e^{\frac{2}{\mu}\frac{\eta}{1-\beta} x} p(x)\mathrm{d}x} \sqrt{\int \left| c(x) - c_n(x) \right|^2 p(x)\mathrm{d}x} + H(u,w,c_n) \\
&= L' \sqrt{\int (x + C')^2 e^{C_X X} p(x)\mathrm{d}x} \left\| c - c_n \right\| + H(u,w,c_n)
\end{aligned} \tag{78}$$

with constants $L', C' > 0$ independent from $n$. Here, to obtain the second line of (78) we used

$$\begin{aligned}
\left| e^{\frac{F(x,u,w,c)}{\mu}} - e^{\frac{F(x,u,w,c_n)}{\mu}} \right| &\le e^{\max\left\{ \frac{F(x,u,w,c)}{\mu}, \frac{F(x,u,w,c_n)}{\mu} \right\}} \left| \frac{F(x,u,w,c)}{\mu} - \frac{F(x,u,w,c_n)}{\mu} \right| \\
&\le e^{\frac{1}{\mu}\left( g(1) + \frac{\lambda U}{\alpha} + \frac{\eta}{1-\beta} X \right)} \left| \frac{F(x,u,w,c)}{\mu} - \frac{F(x,u,w,c_n)}{\mu} \right|
\end{aligned} \tag{79}$$

due to (38). Consequently, we have

$$H(u,w,c) \le C'' \left\| c - c_n \right\| + H(u,w,c_n) \tag{80}$$

with a constant $C'' > 0$ independent from $n$. By virtue of $c_n$, we have $\liminf_{n \to +\infty} \left\| c - c_n \right\| = 0$ and

$$H(u,w,c) \le \liminf_{n \to +\infty} \left\{ C'' \left\| c - c_n \right\| + H(u,w,c_n) \right\} = \liminf_{n \to +\infty} H(u,w,c_n), \tag{81}$$

and hence the lower semi-continuity of $H$ with respect to $c \in \mathfrak{C}$. The lower semi-continuity for $(u,w)$ is proven in an analogous way.

□

***Proof of Proposition 6***

The cases $\kappa = 0,1$ are trivial. If $\kappa \in (0,1)$, then

$$\begin{aligned} G\left(\kappa\theta_1+(1-\kappa)\theta_2\right) &= \ln\left(\int e^{\kappa\theta_1(x)+(1-\kappa)\theta_2(x)}p(x)\mathrm{d}x\right) \\ &\le \ln\left\{\left(\int e^{\theta_1(x)}p(x)\mathrm{d}x\right)^{\kappa}\left(\int e^{\theta_2(x)}p(x)\mathrm{d}x\right)^{1-\kappa}\right\} \\ &= \ln\left(\int e^{\theta_1(x)}p(x)\mathrm{d}x\right)^{\kappa}+\ln\left(\int e^{\theta_2(x)}p(x)\mathrm{d}x\right)^{1-\kappa} \\ &= \kappa\ln\left(\int e^{\theta_1(x)}p(x)\mathrm{d}x\right)+(1-\kappa)\ln\left(\int e^{\theta_2(x)}p(x)\mathrm{d}x\right) \end{aligned}, \quad (82)$$

where the first inequality comes from an application of the classical Hölder's inequality ( $\kappa^{-1}>1$ and $(1-\kappa)^{-1}>1$, and $\kappa+(1-\kappa)=1$):

$$\begin{aligned} \int e^{\kappa\theta_1(x)+(1-\kappa)\theta_2(x)}p(x)\mathrm{d}x &= \int\left(e^{\theta_1(x)}\right)^{\kappa}\left(e^{\theta_2(x)}\right)^{1-\kappa}p(x)\mathrm{d}x \\ &\le \left(\int\left\{\left(e^{\theta_1(x)}\right)^{\kappa}\right\}^{\frac{1}{\kappa}}p(x)\mathrm{d}x\right)^{\kappa}\left(\int\left\{\left(e^{\theta_2(x)}\right)^{1-\kappa}\right\}^{\frac{1}{1-\kappa}}p(x)\mathrm{d}x\right)^{1-\kappa} \\ &= \left(\int e^{\theta_1(x)}p(x)\mathrm{d}x\right)^{\kappa}\left(\int e^{\theta_2(x)}p(x)\mathrm{d}x\right)^{1-\kappa} \end{aligned}. \quad (83)$$

Notice that

$$G(\theta_1)\le G(\theta_2) \text{ if } \theta_1\le\theta_2 \text{ a.s. in } \Omega. \quad (84)$$

Secondly, the mapping $F(x,u,w,c)$ of (32) is convex with respect to $u,w\ge 0, c\in\mathfrak{C}$. Indeed, the three terms are convex, and the sum of finitely many convex functions is again convex. The convexity of $G$ is therefore proven.

Due to the increasing and convex nature of $G$ and the convexity of $F$, for any $\kappa\in(0,1)$ and $u_i,w_i\ge 0, c_i\in\mathfrak{C}$ ($i=1,2$), we obtain

$$\begin{aligned} &H\left(\kappa u_1+(1-\kappa)u_2,\kappa w_1+(1-\kappa)w_2,\kappa c_1+(1-\kappa)c_2\right) \\ &= \mu\ln\mathbb{E}\left[e^{\frac{F\left(X,\kappa u_1+(1-\kappa)u_2,\kappa w_1+(1-\kappa)w_2,\kappa c_1+(1-\kappa)c_2\right)}{\mu}}\right] \\ &\le \mu\ln\mathbb{E}\left[e^{\frac{\kappa F(X,u_1,w_1,c_1)+(1-\kappa)F(X,u_2,w_2,c_2)}{\mu}}\right] \\ &= \kappa\mu\ln\mathbb{E}\left[e^{\frac{\kappa F(X,u_1,w_1,c_1)}{\mu}}\right]+(1-\kappa)\mu\ln\mathbb{E}\left[e^{\frac{F(X,u_2,w_2,c_2)}{\mu}}\right] \\ &= \kappa H(u_1,w_1,c_1)+(1-\kappa)H(u_2,w_2,c_2) \end{aligned}. \quad (85)$$

The convexity of $-\lambda u+\eta w$ is trivial since it is an affine function. Consequently, $-\lambda u+\eta w+H(u,w,c)$ of the problem (37) is convex.

□

***Proof of Proposition 7***

$H_\tau$ is a sum of strictly convex functions and $G$ is convex. Therefore, as for the discussion to obtain (85), for any $\kappa\in(0,1)$ and $u_i,w_i\ge 0, c_i\in\mathfrak{C}$ ($i=1,2$) with $(u_1,w_1,c_1)\ne(u_2,w_2,c_2)$, we obtain

$$
\begin{aligned}
&H_\tau\left(\kappa u_1+(1-\kappa)u_2,\kappa w_1+(1-\kappa)w_2,\kappa c_1+(1-\kappa)c_2\right)\\
&=\mu\ln\mathbb{E}\left[e^{\frac{F_\tau\left(X,\kappa u_1+(1-\kappa)u_2,\kappa w_1+(1-\kappa)w_2,\kappa c_1+(1-\kappa)c_2\right)}{\mu}}\right]\\
&<\mu\ln\mathbb{E}\left[e^{\frac{\kappa F_\tau\left(X,u_1,w_1,c_1\right)+(1-\kappa)F_\tau\left(X,u_2,w_2,c_2\right)}{\mu}}\right]\\
&=\kappa\mu\ln\mathbb{E}\left[e^{\frac{\kappa F_\tau\left(X,u_1,w_1,c_1\right)}{\mu}}\right]+(1-\kappa)\mu\ln\mathbb{E}\left[e^{\frac{F_\tau\left(X,u_2,w_2,c_2\right)}{\mu}}\right]\\
&=\kappa H_\tau\left(u_1,w_1,c_1\right)+(1-\kappa)H_\tau\left(u_2,w_2,c_2\right)
\end{aligned}
\tag{86}
$$

and thus (46). Consequently, the problem (44) is an optimization problem of a strictly convex objective, it has a unique global minimizer (Corollary 2 in Chapter 2 of Zeidler (1995)). The compactness of optimal $u,w$ follows as in **Proposition 4**.

□

***Proof of Proposition 8***

A solution to (31) is denoted as $Z_0^*\in D\times\mathfrak{C}$ with a sufficiently large compact set $D$. The unique solution to the optimization problem (44) is symbolically denoted as $Z_\tau^*\in D\times\mathfrak{C}$ ($\tau>0$). The objectives corresponding to (31) and (44) are denoted as $J_0=J_0(Z)$ and $J_\tau=J_\tau(Z)$, respectively.

For the first statement, for each $(u,w)\in D, c\in\mathfrak{C}$, we obtain

$$
\begin{aligned}
\left|H(u,w,c)-H_\tau(u,w,c)\right|&=\left|\mu\int e^{\frac{F(x,u,w,c)}{\mu}}p(x)\mathrm{d}x-\mu\int e^{\frac{F_\tau(x,u,w,c)}{\mu}}p(x)\mathrm{d}x\right|\\
&\le\mu\int\left|e^{\frac{F(x,u,w,c)}{\mu}}-e^{\frac{F_\tau(x,u,w,c)}{\mu}}\right|p(x)\mathrm{d}x\\
&\le\mu\int e^{\frac{1}{\mu}\left(g(1)+\frac{\lambda U}{\alpha}+\frac{\eta}{1-\beta}x\right)}\left|\frac{F(x,u,w,c)}{\mu}-\frac{F_\tau(x,u,w,c)}{\mu}\right|p(x)\mathrm{d}x\\
&=\int e^{\frac{1}{\mu}\left(g(1)+\frac{\lambda U}{\alpha}+\frac{\eta}{1-\beta}x\right)}\left|F(x,u,w,c)-F_\tau(x,u,w,c)\right|p(x)\mathrm{d}x
\end{aligned}
\tag{87}
$$

By (43), we have the estimate

$$
\begin{aligned}
&\left|F(x,u,w,c)-F_\tau(x,u,w,c)\right|\\
&=\left|\begin{array}{l}g\left(\left|c(x)-\hat{c}\right|\right)+\frac{\lambda}{\alpha}\max\left\{u-\left(1-c(x)\right)x,0\right\}+\frac{\eta}{1-\beta}\max\left\{\left(1-c(x)\right)x-w,0\right\}\\-\left(g\left(\left|c(x)-\hat{c}\right|\right)+\frac{\lambda}{\alpha}m_\tau\left(u-\left(1-c(x)\right)x\right)+\frac{\eta}{1-\beta}m_\tau\left(\left(1-c(x)\right)x-w\right)\right)\end{array}\right|\\
&\le\left(\frac{\lambda}{\alpha}+\frac{\eta}{1-\beta}\right)\tau
\end{aligned}
\tag{88}
$$

whose last line does not depend on $x,u,w,c$. Hence, by **Assumption** 1 we obtain

$$
\left|-\lambda u+\eta w+H(u,w,c)-\left(-\lambda u+\eta w+H_\tau(u,w,c)\right)\right|\le\left(\frac{\lambda}{\alpha}+\frac{\eta}{1-\beta}\right)\tau\int e^{\frac{1}{\mu}\left(g(1)+\frac{\lambda U}{\alpha}+\frac{\eta}{1-\beta}x\right)}p(x)\mathrm{d}x\equiv\bar{C}\tau\,. \tag{89}
$$

Therefore, we have the uniform convergence

$$
\lim_{\tau\to+0}\left|-\lambda u+\eta w+H(u,w,c)-\left(-\lambda u+\eta w+H_\tau(u,w,c)\right)\right|=0\,. \tag{90}
$$

This implies

$$
-\bar{C}\tau\le J_0(Z)-J_\tau(Z)\le\bar{C}\tau\,,\quad\tau>0\,. \tag{91}
$$

From the left inequality, we obtain $-\bar{C}\tau + W_\tau \le W$ , and the from the right one we obtain $W \le W_\tau - \bar{C}\tau$ , yielding $|W - W_\tau| \le \bar{C}\tau$ (this is the last statement of the proposition). Thereby, $W_\tau \to W$ as $\tau \to +0$ .

For the second statement, because each $Z_\tau^*$ is bounded for each $N \in \mathbb{N}$ , and hence has a weak accumulation point in $D \times \mathfrak{D}$, which is denoted as $Z_0^{**}$ by choosing a weakly converging subsequence that is still denoted as $Z_\tau^*$. Because of the optimality of $Z_\tau^*$ and (91), for each $Z \in D \times \mathfrak{C}$, we have

$$J_0(Z^*) = \lim_{\tau \to +0} J_\tau(Z^*) \ge \lim_{\tau \to +0} J_\tau(Z_\tau^*) \ge \liminf_{\tau \to +0} J_\tau(Z_\tau^*) \tag{92}$$

and hence

$$J_0(Z^*) \ge \liminf_{\tau \to +0} J_\tau(Z_\tau^*) . \tag{93}$$

Because of (91), we further obtain

$$J_0(Z_0^{**}) \le \liminf_{\tau \to +0} J_0(Z_\tau^*) \le \liminf_{\tau \to +0} J_\tau(Z_\tau^*) + \liminf_{\tau \to +0} \bar{C}\tau = \liminf_{\tau \to +0} J_\tau(Z_\tau^*) . \tag{94}$$

Combining (93) and (94) yields

$$W = J_0(Z_0^*) \ge J_0(Z_0^{**}) , \tag{95}$$

showing that $W = J_0(Z_0^{**})$ and $Z_0^* = Z_0^{**}$ by the optimality and strict convexity, proving the statement.

□

In addition to these proofs, we prove the following used in the main text.

***Lemma A.1***
*For a non-negative random variable* $X$ *,*

$$\frac{1}{a}\ln\left(\mathbb{E}\left[\exp(aX)\right]\right) \le \frac{1}{b}\ln\left(\mathbb{E}\left[\exp(bX)\right]\right) \quad \textit{for} \quad 0 < a \le b \tag{96}$$

*if* $\mathbb{E}\left[\exp(bX)\right] < +\infty$ .

***Proof of Lemma A.1***
By the classical Jensen's inequality for a convex function $x \ln x$ ( $x > 0$ ), we have

$$\mathbb{E}[Y \ln Y] \ge \mathbb{E}[Y] \ln \mathbb{E}[Y] \tag{97}$$

for a positive random variable $Y$ . Substituting $Y = \exp(aX)$ into (97) yields

$$\mathbb{E}\left[aX \exp(aX)\right] \ge \mathbb{E}\left[\exp(aX)\right] \ln \mathbb{E}\left[\exp(aX)\right] \tag{98}$$

or equivalently

$$-\ln \mathbb{E}\left[\exp(aX)\right] + a\frac{\mathbb{E}\left[X\exp(aX)\right]}{\mathbb{E}\left[\exp(aX)\right]} \ge 0 . \tag{99}$$

Consider a function

$$f(a) = \frac{1}{a}\ln \mathbb{E}\left[\exp(aX)\right] \quad \text{for} \quad a > 0 \tag{100}$$

if $\mathbb{E}\left[\exp(aX)\right] < +\infty$ . By (99), we obtain

$$\begin{aligned} \frac{\mathrm{d}f(a)}{\mathrm{d}a} &= -\frac{1}{a^2}\ln \mathbb{E}\left[\exp(aX)\right] + \frac{1}{a}\frac{\mathbb{E}\left[X\exp(aX)\right]}{\mathbb{E}\left[\exp(aX)\right]} \\ &= \frac{1}{a^2}\left(-\ln \mathbb{E}\left[\exp(aX)\right] + a\frac{\mathbb{E}\left[X\exp(aX)\right]}{\mathbb{E}\left[\exp(aX)\right]}\right) \quad \text{for} \quad a > 0 . \\ &\ge 0 \end{aligned} \tag{101}$$

Then, (96) follows from (101).

□

**Appendix B: Statistics and model identification**
Statistics, which are average, variance, skewness, kurtosis are summarized in **Table B.1**, demonstrating that the statistics are highly non-Gaussian and right-skewed. The streamflow in Case 2 has a higher average and variance than that in Case 1, while their skewness and kurtosis are comparable with each other.

For the model identification of the supCBI process, we used the method of Yoshioka et al. (2023a) where $\pi(\mathrm{d}\rho)$ is identified using a least-square method between empirical and theoretical autocorrelation functions. Another, least-square method is applied to the identification of $v(\mathrm{d}z)$ so that the sum of the squares of the relative errors of the four statistics in **Table B.1** is minimized. In this paper, a tempered stable case $v(\mathrm{d}z)=z^{-(\alpha_v+1)}\exp(-\beta_v z)\mathrm{d}z$ ( $\alpha_v<1$ , $\beta_v>0$ ) and the Gamma case $\pi(\mathrm{d}\rho)\sim z^{\alpha_\pi-1}\exp(-\rho/\beta_\pi)\mathrm{d}\rho$ ( $\alpha_\pi>1$ , $\beta_\pi>0$ ) were assumed. Advantages of these parameterizations are that the mathematical justification discussed in the main text are applicable and the autocorrelation function of a sub-exponential form is found as $\left(1/\left(1+\beta_\pi\left(1-BM_1\right)\tau\right)\right)^{\alpha_\pi-1}$ with $\tau\geq 0$ being the time lag. **Table B.2** presents identified parameter values for Cases 1 and 2. The empirical and theoretical autocorrelation functions agree well as demonstrated in **Figures B.1-B.2** below, and the PDFs have also been fitted well as demonstrated in the main text. Interestingly, the streamflow in Case 1 is of the long-memory type because the integral of the autocorrelation function diverges near $\tau=0$ due to $\alpha_\pi\leq 2$, while that in Case 2 does not due to $\alpha_\pi>2$. Nevertheless, both cases exhibit $\alpha_\pi$ close to 2 as in the data of the other river environment (Yoshioka, 2022), suggesting some physical universality that is beyond the scope of this paper but would be of theoretical interest.

**Table B.1.** Statistics at the study sites of Cases 1 and 2.

| Statistics | Case 1 | Case 2 |
|---|---|---|
| Average ($\mathrm{m^3/s}$) | 12.47 | 37.17 |
| Variance ($\mathrm{m^6/s^2}$) | 340.57 | 3547.31 |
| Skewness (-) | 10.91 | 11.85 |
| Kurtosis (-) | 209.23 | 218.0 |

**Table B.2.** Statistics at the study sites of Cases 1 and 2.

| Parameter | Case 1 | Case 2 |
|---|---|---|
| $\underline{X}$ ($\mathrm{m^3/s}$) | 0.00 | 5.86 |
| $\alpha_\pi$ (-) | 1.82 | 2.09 |
| $\beta_\pi$ (1/hour) | 0.0686 | 0.0783 |
| $A$ ( $\mathrm{hour}\cdot\mathrm{s}^{-\alpha_v}\cdot\mathrm{m}^{3\alpha_v}$ ) | 0.0300 | 0.0937 |
| $B$ ( $\mathrm{s}^{1-\alpha_v}\cdot\mathrm{m}^{-3(1-\alpha_v)}$ ) | 0.0285 | 0.0236 |
| $\alpha_v$ (-) | 0.852 | 0.803 |
| $\beta_v$ ($\mathrm{s/m^3}$) | 0.00450 | 0.00144 |

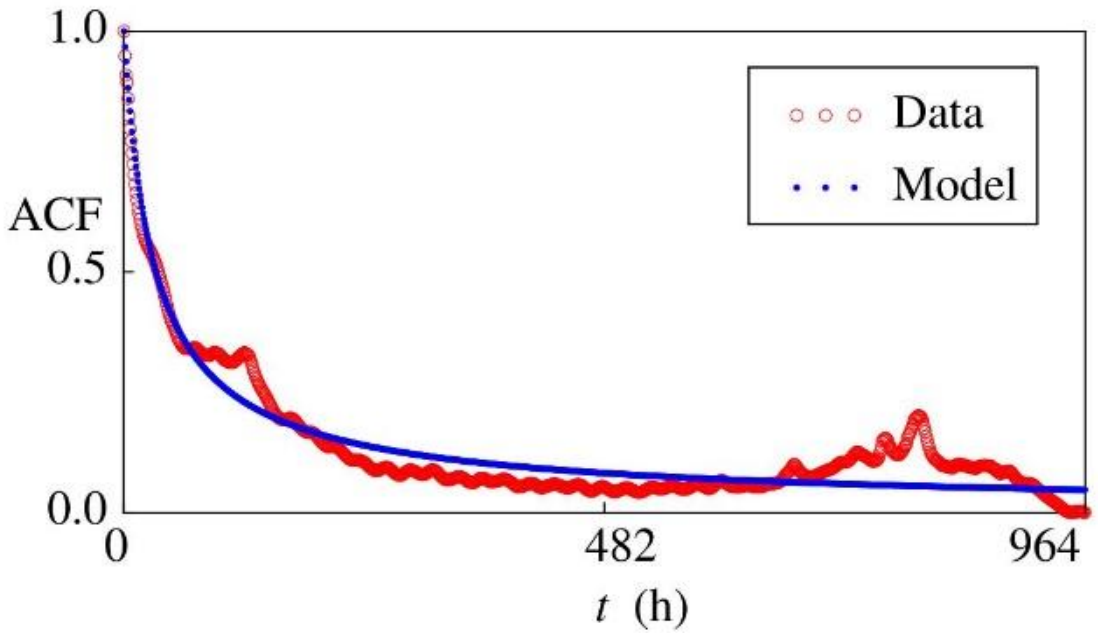

**Figures B.1.** Comparison of empirical and theoretical autocorrelation functions (ACFs) in Case 1.

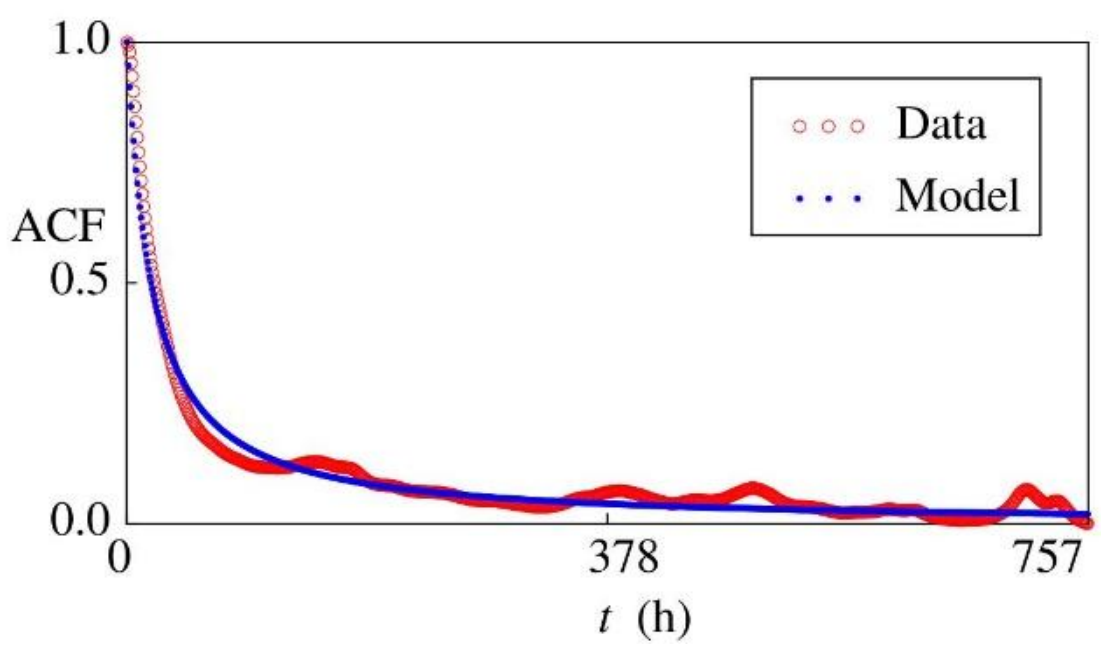

**Figures B.2.** Comparison of empirical and theoretical autocorrelation functions (ACFs) in Case 2.

**Appendix C Monte-Carlo method**

The algorithm of the Monte-Carlo method used for simulating sample paths of supCBI processes is explained in this appendix. It is **Algorithm 2** based on the Markovian lift. We simulate sample paths of mutually-independent CBI processes having specified reversion speed $\rho_i$. This algorithm is a straightforward extension of the algorithm for supOU processes (Yoshioka et al., 2022b). After simulating sample paths of the supCBI processes, we sum up them with the weights $c_i$. In **Step 2**, the algorithm of Algorithm 0 in Kawai and Masuda (2011) with a common Euler-Maruyama discretization is employed to generate tempered stable jump. The time step $\Delta t$ in the simulation is 0.01/24 (day).

**Algorithm 2**

***Step 0*** *Prepare sequences* $\{c_i\}_{1\le i\le n}$ *and* $\{\rho_i\}_{1\le i\le n}$*, the discrete time step* $\Delta t > 0$. *Set* $i = 1$.

***Step 1*** *Simulate the discrete-time approximation of the CBI process with the parameters* $(c_i, \rho_i)$ *for time steps* $j = 0,1,2,...$ *with an initial condition* $Y_0(\rho_i) \ge 0$*, and obtain* $\{Y_{j\Delta t}(\rho_i)\}_{j=0,1,2,...}$.

***Step 2*** *For each* $j = 0,1,2,...$*, compute* $\left\{\underline{X} + \sum_{i=1}^{n} Y_{j\Delta t}(\rho_i)\right\}_{j=0,1,2,...}$.